\documentclass[12pt]{amsart} \usepackage{hyperref,amscd}

\usepackage{amsmath,amssymb,epsfig,amscd}

\newcommand{\new}{\newcommand}

\new{\abs}[1]{\left|#1 \right|}
\new{\norm}[1]{\left\| #1 \right\|} 
\new{\bracket}[1]{\langle #1
\rangle} 
\new{\defequals}{\stackrel{\operatorname{def}}{=}}
\new{\into}{Hookrightarrow} 
\new{\comb}[2]{{#1 \choose #2}}
\new{\tensor}{\otimes}
\new{\tnsr}{\otimes} 
\new{\iso}{\cong} 
\new{\union}{\cup}
\new{\maps}{\colon} 
\new{\goesto}{\rightarrow}
\new{\R}{{\mathbb{R}}} 
\new{\C}{{\mathbb{C}}}

\new{\pic}[5]{\raisebox{#3pt}{
Hspace{#4pt}\epsfig{file=#1.eps,height=#2pt}Hspace{#5pt}}}

\newtheorem{proposition}{Proposition}
\newtheorem{theorem}{Theorem}

\newtheorem{lemma}{Lemma} 

\newtheorem{corollary}{Corollary} 

\newtheorem{definition}{Definition}

\newtheorem{remark}{Remark}

\newcounter{letter}
\newenvironment{alist}{
\begin{list}{(\alph{letter})}{\usecounter{letter}}
}{\end{list}}

\newenvironment{pf}{
\begin{trivlist} \item[Proof:]}{\qed \end{trivlist}}

\new{\lieg}{\mathfrak{g}}
\new{\lieh}{\mathfrak{h}}
\new{\hath}{\Check{h}}
\new{\NN}{{\mathbb{N}}}
\new{\ZZ}{{\mathbb{Z}}}
\new{\QQ}{{\mathbb{Q}}}
\new{\CC}{{\mathbb{C}}}
\new{\trunc}{\hat{\tensor}}
\new{\Weyl}{{\mathcal W}}
\new{\A}{{\mathcal A}}
\new{\Cat}{{\mathcal C}}
\new{\s}{{\mathbf s}}
\new{\q}{{\mathbf q}}
\new{\D}{{\mathcal D}}
\new{\HH}{{\mathcal H}}
\new{\WW}{{\mathfrak{W}}}
\new{\semidirect}{\ltimes}

\new{\tr}{\operatorname{tr}} 
\new{\qdim}{\operatorname{qdim}} 
\new{\qtr}{\operatorname{qtr}} 
\new{\range}{\operatorname{Range}} 
\new{\dom}{\operatorname{Domain}}
\new{\Lie}{\operatorname{Lie}} 
\new{\Sym}{\operatorname{Sym}} 
\new{\Fun}{\operatorname{Fun}} 
\new{\Aut}{\operatorname{Aut}}
\new{\Ad}{\operatorname{Ad}}
\new{\ad}{\operatorname{ad}}
\new{\coad}{\operatorname{coad}}
\new{\vol}{\operatorname{vol}}
\new{\Int}{\operatorname{Int}}
\new{\Ext}{\operatorname{Ext}}
\new{\Image}{\operatorname{Im}}
\new{\Ker}{\operatorname{Ker}}
\new{\Map}{\operatorname{Map}}
\new{\antisym}{\operatorname{Antisym}}
\new{\Hom}{\operatorname{Hom}}
\new{\Null}{\operatorname{Null}}
\new{\End}{\operatorname{End}}
\new{\Stab}{\operatorname{Stab}}

\begin{document}

\title[Roots of Unity]{Quantum Groups at Roots of Unity and Modularity}

\author{Stephen F. Sawin}
        \address{ Department of Math and C. S. \\ Fairfield University \\ Fairfield, CT
06824-5195}
    \email{ sawin@cs.fairfield.edu} 
 
\begin{abstract}
We develop the basic
representation theory of all quantum groups at all roots of unity, including
Harish-Chandra's Theorem, which allows us to show that an appropriate
quotient of a subcategory gives a semisimple ribbon category. This
work generalizes previous work on the 
foundations of representation theory of quantum groups at roots of
unity which applied
only to quantizations of the simplest groups, or to certain fractional
levels, or only to the projective form of the group. The second half
of this paper applies the representation to give a
sequence of results crucial to applications in topology.  In
particular  for each
compact, simple, simply-connected Lie group we show that at each integer 
level the quotient categopry is in fact modular (thus leading to a
Topological Quantum Field Theory),  we determine
when at fractional levels the corresponding category is modular, and  we
 give a quantum version of the Racah formula for the decomposition
of the tensor product.  
\end{abstract}

 \maketitle    
     
\tableofcontents

\section*{Introduction}
In \cite{witten89a} Witten argued that Chern-Simons theory for a
compact, connected, simply-connected simple Lie group at integer
level $k$ should
yield an invariant of links in a (biframed) three-manifold.  He also
sketched how 
 to compute this
invariant combinatorially  using two-dimensional conformal field
theory, and worked out the $\mathrm{SU}(2)$ invariant in enough 
detail   to demonstrate
that if  well-defined it would have to give the Jones
polynomial.  in \cite{RT91}, Reshetikhin
and Turaev constructed an invariant that met all of Witten's
criteria using the quantum group associated to $\mathfrak{sl}_2$ (the
complexification of the Lie algebra of  $\mathrm{SU}(2)$)  with
the quantum parameter equal to a root of unity depending on the level $k.$ At
this point an overall program was clear.  To each simple Lie algebra
there was associated a quantum group.  By understanding the
representation theory of this quantum group at roots of unity in an
analogous fashion to Reshetikhin and Turaev's work on
$\mathfrak{sl}_2,$ one could presumably show that this representation theory
formed a \emph{modular tensor category,} and thus construct an
invariant of links and three-manifolds, presumably the one
Witten associated to the corresponding compact, simple Lie group.  

Remarkably, the intervening twelve years have not sufficed to complete
this apparently straightforward program.  The difficulty is even more
surprising considering that around the same time, an analogy was
identified between quantum groups at roots of unity and modular
groups, which brought both the attention of many skilled algebraists
and a host of useful techniques to bear on questions of the
representation theory of quantum groups at roots of unity.

Much of the obstacle to the complete resolution of this problem
appears to be faulty communication between those working on the algebraic
questions and the topologists and mathematical physicists interested
in the link and three-manifold invariants.  There are two basic
confusions.  First, most of the algebraic work has been done at 
roots of unity of odd degree (in fact, degree prime to the entries of
the Cartan matrix), 
while the values that correspond to integer levels and 
hence the cases of primary interest for their relation to physics are
 roots of unity of even degree (in fact degrees which are  multiples
 of the entries of the Cartan 
matrix).  Second, much of the algebraic work deals only with
representations with highest weights in the root lattice, while in
topology and physics one is interested in all representations whose
highest weights lie in the weight lattice (indeed the most
straightforward invariants geometrically, those coming from the
simply-connected groups, require  consideration of all the
representations, though most of the information is available from the
invariants depending only on root lattice representations
\cite{Sawin02a,Sawin02b}).  Thus remarkably, there is no proof in the
literature that from any quantum group at roots of unity one can
construct three-manifold invariants, or even semisimple ribbon
categories, though of course these facts are widely understood to be true.

To explain the first confusion requires a careful understanding of the
parameter in the definition of the quantum group.  Since this
parameter always appears with an  exponent which is the
inner product of elements of the dual to the Cartan subalgebra, it
amounts to a choice of scale for this inner product.  We work with a
parameter $q$ normalized so that the definition of the quantum group
and the basic representation theory require only integral powers of
$q.$  Our $q$ corresponds to that of most authors, including Lusztig
\cite{Lusztig88,Lusztig89,Lusztig90b,
Lusztig90c,Lusztig90d,Lusztig92c,Lusztig93}    
(who sometimes calls it $v$), Andersen et al \cite{APWX91,
AW92,APW92,Andersen92,Andersen93,Andersen94,AJS94,
Andersen95,APW95,AP95,Andersen97a,Andersen97b,Andersen00,AP00,Andersen01},
Kirillov \cite{Kirillov96} and Chari and Pressley \cite{CP94}.
However, because ${\mathfrak sl}_2$ uniquely among Lie algebras has all
entries in its Cartan matrix divisible by $2,$ quantum
$\mathfrak{sl}_2$ can be defined with only integral powers of $q^2,$
and thus it 
is conventional when considering only  $\mathfrak{sl}_2$ to refer to
what we call $q^2$ as $q.$  This is the convention 
in, for example, Reshetikhin and Turaev \cite{RT90,RT91}, Kirby and
Melvin \cite{KM91}. Our $q^2$ is also the $t$ of the Jones polynomial
\cite{Jones85} and
hence $A^4$ in Kauffman's bracket \cite{Kauffman87a,Kauffman90a}.  In our
normalization the level $k$ 
of Witten's invariant corresponds to quantum groups at level
$e^{\pi i/D(k + \check{h})}$ where $\check{h}$ is the dual Coxeter
number, and $D$ is the ratio of the square of the length of a long
root to that of a short root (i.e., the biggest absolute value of an
off-diagonal entry of the Cartan matrix).   In
particular it is an $l$th root of unity where $l$ is a multiple of
$2D.$ 

Further confusing the issue is that when the determinant of the Cartan
matrix is not  $1,$ the $R$-matrix cannot be defined
(at least in the presence of the weight lattice representations)
without introducing fractional powers of $q.$  We recover integral
powers if we express everything in terms of the parameter $s,$ where
$q=s^L,$ with $L$ the index of the root lattice in the weight
lattice.  Also,
some authors use a slightly different set of 
generators which gives a more symmetric presentation but requires
half-integer powers of $q$ and $s,$ including Rosso, Kirby and Melvin,
and Kirillov 
\cite{Rosso88,Rosso90,KM91,Kirillov96}.

On the algebraic side, the basic facts of quantum groups at roots of
unity are somewhat easier if one assumes that the degree of the root of
unity is odd and prime to $D.$ Much of the work in this field has
focused on the relationship between the representation theory for
algebraic groups over a 
field of prime characteristic $p$ and the
representation of the corresponding quantum group at a $p$th root of
unity.  Consequently the restriction on the degree of the root of unity
seemed harmless and most early  work in the field maintained that
restriction \cite{Lusztig89,APWX91,
AW92,APW92,Andersen92,Andersen93,Andersen94,AJS94,
Andersen95,APW95,CP94,Paradowski94}.  Thus the results in these works
apply to an 
entirely disjoint set of situations from the quantum groups that arise
in connection with affine Lie algebras, conformal field theory and
three-manifold invariants. 

In fact, for the elementary sorts of results about the representation
theory of quantum groups that are necessary to construct the
three-manifold invariants and similar tasks, the restrictions on the
degree of the root of unity do not appear to be essential, and in fact
Andersen and Paradowski \cite{AP95}  reproduce many of
these results in the general situation, in particular the quantum
version of the Racah formula for the tensor product of
representations.  However, Andersen and Paradowski restrict attention to
representations whose weights lie in the root lattice.  This restriction
corresponds
on the Lie group level to considering only the adjoint form of the
group.  To understand the geometry of the construction it would be
best to be able to construct the invariant for all forms of the group.

This article redresses this lack, proving the fundamental results on
the representation of quantum groups needed for applications to
topology and physics for all nongeneric values of the parameter and all
representations.  The proofs rely very little on technology borrowed
from algebraic groups over finite fields, and  as such one can hope
they will
be more accessible to researchers interested in quantum
invariants.  The tools are those parts of classical Lie algebra theory
that generalize directly to quantum groups, up to the PBW theorem and
 Harish-Chandra's theorem, together with the
$R$-matrix and Lusztig's integral form. 

Section 1 fixes notation and defines the general quantum group
$U_q(\lieg)$ over $\QQ(q),$ reviewing some elementary results from
the literature.  The key nontrivial results we will use are the quantum
version of the 
Poincar\'{e}-Birkhoff-Witt Theorem and the existence of the integral
form $U^{\text{res}}_\A(\lieg),$  both due to Lusztig. Also in the
section is the construction of the integral form
$U^\dagger_{\A'}(\lieg),$ which is literally a ribbon Hopf algebra.  The
existence of such a form is new, and although the construction is not
deep, this form is crucial at several points in the article, and seems
interesting in its own right.

Section 2 constructs the quantum group at roots of unity (two forms,
the standard $U_\q^{\text{res}}(\lieg)$ and the ribbon Hopf algebra
$U_\s^\dagger(\lieg)$).  It also defines the affine Weyl group and
uses it to prove Theorem~\ref{harish_chandra}, the quantum version of
Harish-Chandra's Theorem (that is, that Weyl modules have the same
character if and only if they are related by the affine Weyl group.
We do not characterize the center of the quantum group, which is also
often referred to as Harish-Chandra's Theorem).  The description
of the affine Weyl group, which is entirely elementary, has one point
worthy of comment, though it is not new.   The action of the affine
Weyl group depends 
subtly on the divisibility of the degree of the root of unity, as
was observed in \cite{AP95}.  At roots of unity corresponding to
integer levels it is what is traditionally called the affine Weyl
group of $\lieg$ in Lie algebra theory and loop groups, while at
certain fractional levels it is the affine Weyl group of the dual Lie
algebra.   The existence of
$(2+1)$-dimensional theories at fractional level with a nonstandard Weyl
alcove has not been widely recognized and does not seem to have a
counterpart in the physics literature.  

Harish-Chandra's Theorem is proven for $U_\s^\dagger(\lieg).$ It is
not in general true for $U_\q^{\text{res}}(\lieg),$ though a  key
consequence, the Linkage Principle, is. The proof follows the
classical proof, except that it relies explicitly on the $R$-matrix,
and at fractional levels supplements the Harish-Chandra map to the
center (which 
is not onto) with additional central elements not in the range.  The
same argument gives Harish-Chandra's theorem (with the ordinary Weyl
group) for $U_{\A'}^\dagger.$

Section 3 (working over $U_\q^{\text{res}}(\lieg)$ and
$U_\s^\dagger(\lieg)$),   defines tilting modules as in
\cite{Paradowski94,AP95}, proves they form a tensor category, and
shows that each is a sum of 
 highest weight tilting modules.
This section expands and simplifies the work in these two references,
using only the elementary theory of quantum groups themselves, rather
than techniques from algebraic groups. 

Section 4 identitifes the so-called negligible modules, which have the
property that all intertwiners from the module to itself have quantum
trace zero.  This implies in particular that the link invariant is
trivial whenever a component is labeled by a negligible module.  It is
shown here that in fact every highest weight tilting module outside
the Weyl alcove is negligible. This section follows the arguments in
\cite{Paradowski94,AP95} closely.

Section 5 gives an application of the technology developed in the
previous sections.  Specifically it gives a quantum version of the
Racah formula (also called the Racah-Speiser formula, see MacFarlane
et al.\cite{MOR67})
which expresses multiplicities of the tensor products of two Weyl
modules in terms of weight multiplicities.  A version of this formula
for Weyl modules whose highest weights are required to be in the root
lattice appears in 
\cite{AP95}. 

Finally Section 6 addresses the link and three-manifold invariants and TQFTs
constructable from $U_\s^\dagger.$  The ribbon category
associated to the set of all modules and the one associated to the set of all tilting modules are
given.  The quotient of the tilting modules by the negligible tilting
modules is given and proven to be a semisimple ribbon category.
In the physically interesting case, namely where $2D$ divides $l,$ the
resulting category is shown to be modular, using the previous sections
and well-known results.  The behavior at other roots of unity is
analyzed and the cases when the theory is modular and when it admits a
TQFT or Spin TQFT are identified.
These results are all new, although the results on the physically
interesting case have been widely understood to be true for years.
 
The details of many of the following constructions will depend heavily
on several parameters depending on the root system, including the ratio $D$ of the
square lengths of the long and short roots, the smallest integer $L$
such that $L$ times any inner product of weights is an integer ($L$
always divides the index of the 
fundamental group, that is of the weight lattice modulo the root
lattice), the Coxeter number $h$ and
the dual Coxeter number $\hath.$  For
convenience we summarize these quantities
(all information taken from \cite{Humphreys72}).

\vspace{6pt}

\begin{tabular}{|r||c|c|c|c|c|c|c|c|c|c|c||}
\hline
& $A_n$ & $B_{2n+1}$ & $B_{2n}$ & $C_n$ & $D_{2n}$ & $D_{2n+1}$ & $E_6$ & $E_7$ &
$E_8$ & $F_4$ & 
$G_2$\\
\hline
$L$& $n+1$ & $2$ & $1$ & $1$ & $2$ & $4$ & $3$ & $2$ & $1$ & $1$ & $1$ \\
\hline
$D$ & $1$ & $2$ & $2$ & $2$ & $1$ &  $1$ & $1$ & $1$ & $1$ & $2$ & $3$ \\
\hline
$h$ &$n+1$ & $4n+2$ & $4n$ & $2n$ & $4n-2$ & $4n$ & $12$ & $18$ & $30$ &
$12$ & $6$ \\
\hline
$\strut{\hath}$ &$n+1$ & $4n+1$ & $4n-1$ & $n+1$ & $4n-2$ & $4n$ & $12$ & $18$ & $30$ &
$9$ & $4$ \\
\hline
\end{tabular}

\vspace{6pt}

The author would like to thank Anna Beliakova for pointing out his
misinterpretation of the results in \cite{AP95}, which led to this
article. 

\section{Quantum Groups}
We follow Humphreys \cite{Humphreys72} for all our results on Lie
algebras, and for the most part, notation.  The following paragraph
follows \cite{Humphreys72} Chapter 10.

Let $\lieg$ be a complex simple Lie algebra, let $\lieh$ be a Cartan
subalgebra and $\lieh^*$ its dual vector space.  Let $\Phi \subset
\lieh^*$ be the root system of $\lieg.$ Let $\bracket{\,\cdot\,,\,
\cdot \,}$ be the unique inner product on $\lieh^*$ (and hence on
$\lieh$) such that $\bracket{\alpha,\alpha}=2$ for every short root
$\alpha \in \Phi $ (this convention guarantees that the inner product
of two roots is an integer.  It differs from $(\,\cdot\,,\, \cdot \,)$
normalized in the convention of physics and loop groups, so that the
longest root has length $2.$ $(\,\cdot\,,\, \cdot \,)=
\bracket{\,\cdot\,,\, \cdot \,}/D,$ where $D=3$ for $G_2,$ $D=2$ for
$B_2,$ $C_n,$ and $F_4$ and $D=1$ otherwise.). Let $\check{\Phi}=\{
\check{\alpha}=2\alpha/\bracket{\alpha,\alpha} \,|\, \alpha \in
\Phi\}$ be the dual root system to $\Phi.$ Let $\Lambda=\{\lambda \in
\lieh^*|\bracket{\lambda,\check{\alpha}} \in \ZZ, \, \forall \alpha
\in \Phi\}$ be the weight lattice,  $\Lambda_r=\ZZ\Phi \subset
\Lambda$ be the root lattice, and
$\check{\Lambda}_r=\ZZ \check{\Phi} \subset \frac{1}{D}\Lambda$ be the
dual root lattice.  Let $
\Weyl,$ the \emph{Weyl group,} be the group of 
isometries of $\lieh^*$ generated by reflections about hyperplanes
perpendicular to the roots
$\alpha \in \Phi.$ Thus in particular for each root $\alpha \in \Phi$
we have a reflection $\sigma_\alpha \in \Weyl$ defined by
$\sigma_\alpha(\lambda)=\lambda-\bracket{\lambda, \check{\alpha}}\alpha.$
In fact we will most often be interested in the translated action of
the Weyl group, which is defined by $\sigma \cdot \lambda=
\sigma(\lambda+\rho)-\rho,$ where $\rho=\sum_{\alpha > 0} \alpha/2.$

Let $L$ be the least integer  such
that $L\bracket{\lambda,\gamma} \in \ZZ$ whenever $\lambda,\gamma \in
\Lambda.$ Let $\Delta=\{\alpha_1, \ldots, \alpha_N\} \subset \Phi$ be
a base, $(a_{ij})=\bracket{\alpha_i,\check{\alpha_j}}$ be the Cartan
matrix, and let $\alpha > \beta$ mean $\alpha-\beta$ is a nonnegative
linear combination of the elements of $\Delta.$ Let
$\Lambda^+=\{\lambda \in \Lambda | \bracket{\lambda,\alpha_i} \geq 0,
\forall \alpha_i \in \Delta\}$ be the set of nonnegative integral
weights.  Let $\theta$ be the longest root, 
i.e., the unique long root in $\Phi \cap \Lambda^+,$ and let $\phi$ be
the unique short root in the same intersection.  Finally, let

$\hath=\bracket{\rho,\check{\theta}}+1$ be the \emph{dual Coxeter
number} and $h=\bracket{\rho,\check{\phi}}+1$ be the \emph{Coxeter number.}

Let $\A=\ZZ[q,q^{-1}].$ Given integers $m,n$ let
$$[n]_q=(q^n-q^{-n})/(q-q^{-1}) \in \A,$$
$$[n]_q!=[n]_q\cdot [n-1]_q
\cdots [1]_q \in \A,$$
$$\left[\begin{matrix} m\\ n \end{matrix}
\right]_{q}=[m]_q!/([n]_q![m-n]_q!) \in \A.$$
 Let
$d_i=\bracket{\alpha_i,\alpha_i}/2$ (so $d_i=1$ for short roots and
$d_i=D$ for long) and let $q_i=q^{d_i}.$ Following
\cite{CP94} (Except we use $E$ and $F$ for their $X^\pm$) we define
the Hopf algebra $U_q(\lieg)$ over $\QQ(q)$ with generators $E_i,$
$F_i,$ $K_i$ (for $1 \leq i \leq N$)
$$K_iK_j=K_jK_i \qquad K_iK_i^{-1}= K_i^{-1}K_i=1,$$
$$K_i E_j K_i^{-1}=q^{\bracket{\alpha_i,\alpha_j}} E_j, \qquad K_i F_j K_i^{-1}=q^{- \bracket{\alpha_i,\alpha_j}} F_j,$$
$$E_iF_j - F_j E_i=\delta_{i,j} \frac{K_i-K_i^{-1}}{q_i-q_i^{-1}},$$
$$\sum_{r=0}^{1-a_{ij}} (-1)^r \left[\begin{matrix} 1-a_{ij}\\ r \end{matrix} \right]_{q_i}(E_i)^{1-a_{ij}-r} E_j (E_i)^r=0 \qquad \text{if } i \neq j,$$
$$\sum_{r=0}^{1-a_{ij}} (-1)^r \left[\begin{matrix} 1-a_{ij}\\ r \end{matrix} \right]_{q_i}(F_i)^{1-a_{ij}-r} F_j (F_i)^r=0 \qquad \text{if } i \neq j,$$
$$\Delta(K_i)=K_i \tensor K_i,$$
$$\Delta(E_i)=E_i \tensor K_i + 1 \tensor E_i, \qquad \Delta(F_i)=F_i \tensor 1 + K_i^{-1} \tensor F_i,$$
$$S(K_i)=K_i^{-1}, \quad S(E_i)=-E_i K_i^{-1}, \quad S(F_i)=-K_i F_i,
$$
$$\epsilon(K_i)=1, \qquad \epsilon(E_i)=\epsilon(F_i)=0.$$

Define $E_i^{(l)}=E_i^l/[l]_{q_i}!$ for $l \in \NN,$ and likewise for
$F_i.$ We define 
the $\A$-subalgebra $U_\A^{\text{res}}(\lieg)$ of $U_q(\lieg)$ to be
generated by the elements $E_i^{(r)},$ $F_i^{(r)},$ $K_i^{\pm 1},$ for
$1 \leq i \leq n$ and $r \geq 1.$ We summarize from \cite{CP94}[9.3,
10.1] some important facts about $U_q(\lieg),$
$U_\A^{\text{res}}(\lieg),$ and their respective representation
theories.

$U_\A^{\text{res}}(\lieg)$ is an integral form of $U_q(\lieg)$ in the
sense that $U_q(\lieg)=U_\A^{\text{res}}(\lieg) \tensor_\A \QQ(q),$
and $U_\A^{\text{res}}(\lieg)$ is a free $\A$-algebra. There exist
$E_{\beta_1}, F_{\beta_1}, E_{\beta_2}, F_{\beta_2}, \ldots ,
E_{\beta_N}, F_{\beta_N}, \in U_\A^{\text{res}}(\lieg),$ where
$\beta_1, \ldots, \beta_N$ is an enumeration of the positive roots,
such that each $E_{\beta_i},F_{\beta_i}$ satisfies $K_j E_{\beta_i} K_j^{-1}=
q^{\bracket{\alpha_j,\beta_i}}E_{\beta_i}$ and the set of all
$(E_{\beta_N})^{(l_N)}\cdots (E_{\beta_1})^{(l_1)}$ forms a basis for
the subalgebra $U^+_q(\lieg)$ generated by $\{1,E_i^{(k)}\}$ (and likewise
for $F,$ with $U^-_q(\lieg)$ defined correspondingly).

 For each $\lambda \in \Lambda^+$ there is a unique (up to
 isomorphism) irreducible,
  $U_q(\lieg)$-module with a vector  $v$ such that
  $K_iv=q^{\bracket{\alpha_i,\lambda}}v$ and no vector $w$ satisfies
  $K_iw=q^{\bracket{\alpha_i,\gamma}}w$ for $\gamma >\lambda.$  This
  module, called the Weyl
 module $W^\lambda_q,$ 
  is a direct sum of its weight spaces and the dimensions of
 these weight spaces are the same as that of the classical Weyl
 module $W^\lambda.$   The tensor product of two Weyl modules is isomorphic to a
 direct sum of Weyl modules with
 multiplicities the same as those in the classical case.  If $v$ is a
 highest weight vector of a module $W^\lambda$ (that is a vector of
 weight $\lambda$)  then
 $W_{\A}^\lambda=U_\A^{\text{res}}(\lieg)\cdot v$ is a
 $U_\A^{\text{res}}(\lieg)$-submodule of $W^\lambda,$  an $\A$-form
 of $W^\lambda,$ and  a direct sum of its intersections with the
 weight spaces of $W^\lambda,$ each of which is a free $\A$-module of
 finite rank.

One of the most important reasons for considering quantum groups is
that they are ribbon Hopf algebras, and give invariants of
links.  It is thus a source of some embarrassment that none of the
versions of the quantum group defined here are ribbon Hopf algebras!
There are three obstacles to defining a ribbon structure
(specifically, to giving an $R$-matrix).  The first was discussed in
the introduction, that to define the $R$-matrix, requires fractional
powers of $q$ (Though in the end the properly normalized link
invariant contains only integer powers of $q,$ see Le \cite{Le00a}). The 
second is that an insufficient subset of the enveloping algebra of the
Cartan subalgebra has been included in these definitions to write down the $R$-matrix, even
formally.  The third is that the $R$-matrix involves an infinite sum,
though in any particular finite-dimensional representation only
finitely many terms are nonzero.  Below we give an integral form of
the quantum group which is a ribbon Hopf algebra (technically a
topological Hopf algebra, in the sense that the comultiplication maps
to a completed tensor product).  

Let $\A'=\ZZ[s,s^{-1}]$ and define a monomorphism $\A \to \A'$ by
$q \mapsto s^L$ (henceforth we will treat this monomorphism as an
inclusion and simply write $q=s^L$).  Define the $\A'$ Hopf algebra
$U_{\A'}^{\text{res}}(\lieg)=U_\A^{\text{res}}(\lieg) \tensor_\A \A'$
and the $U_{\A'}^{\text{res}}(\lieg)$-module
$W_{\A'}^\lambda=W_{\A}^\lambda \tensor_\A \A'.$ 

Now consider the Hopf algebra of functions on the additive group
$\Lambda.$ The collection of all set-theoretic functions from
$\Lambda$ to $\A',$ $\Map(\Lambda,\A'),$ is naturally an algebra over
$\A'$ with
pointwise multiplication.  It is a (topological) Hopf algebra when
given the comultiplication $\Delta(f)(\mu,\mu')= f(\mu+\mu'),$ the
counit $\epsilon(f)=f(0),$ and the antipode $S(f)(\mu)=f(-\mu)$ for $f
\in \Map(\Lambda,\A')$ and $\mu,\mu' \in \Lambda.$  Here $$\Delta
\colon \Map(\Lambda, \A') \to \Map(\Lambda \times \Lambda, \A').$$  The
latter space contains the natural embedding of $\Map(\Lambda,\A')
\tensor \Map(\Lambda,\A')$ as a dense subspace in the topology of
pointwise convergence, and thus may be viewed as the completed tensor
product.  A topological basis for this Hopf algebra of functions is given by
$\{\delta_\lambda\}_{\lambda \in \Lambda},$ where
$\delta_\lambda(\gamma)=\delta_{\lambda,\gamma}.$  By topological
basis we mean the elements are linearly independent and span a dense subspace
of $\Map(\Lambda,\A')$ in the topology of pointwise convergence.

Recall any Abelian group with a homomorphism to its dual has an
$R$-matrix associated to the homomorphism in the  Hopf algebra of functions.
In the case of the homomorphism $\lambda \mapsto
s^{L\bracket{\lambda,\,\cdot\,}},$ the 
$R$-matrix is $\sum_{\lambda,\gamma} s^{L\bracket{\lambda,\gamma}} 
\delta_\lambda \tensor \delta_\gamma,$ which once again is an element
not of the tensor product of the Hopf algebra with itself but of the
completion. Defining $\lambda_i \in \Lambda$ by
$\bracket{\lambda_i,\check{\alpha_j}}=\delta_{i,j}$ we can write the
canonical dual element to the pairing as $\sum_i \lambda_i \tensor
\check{\alpha_i},$ and then being somewhat abusive of notation can
refer to the $R$-matrix above 
as
$q^{\sum_i \lambda_i \tensor
\check{\alpha_i}}.$  In the same vein we shall use $q^\lambda$ to
refer to the homomorphism $\sum_{\gamma \in \Lambda}
s^{L\bracket{\lambda,\gamma}} \delta_\gamma.$ Let
$U^\dagger_{\A'}(\lieh)$ be
$\Map(\Lambda, \A')$ viewed as a topological
ribbon Hopf algebra.

$U_{\A'}^\dagger(\lieh)$ acts on $U_{\A'}^{\text{res}}(\lieg)$ via the
$\Lambda$-grading of 
$U_{\A'}^{\text{res}}(\lieg).$  Specifically,  define the weight
of a monomial in $\{E_i,F_i,K_i\}$ to be the sum of $\alpha_i$ for
each factor of $E_i$ and $-\alpha_i$ for each factor of $F_i.$ Then $f
\in U_{\A'}^\dagger(\lieh)$ acts on a monomial $X$ by
$f[X]=f(\text{weight}(X))X$ and extends linearly.  This action is 
automorphic in the sense that  $\Delta(h[X])=h^{(1)}[X^{(1)}] \tensor
h^{(2)}[X^{(2)}]$ and $h[XY]=h^{(1)}[X]h^{(2)}[Y]$ (here we use
Sweedler's notation, writing $\Delta(h)$ as $\sum_\beta
h^{(1)}_\beta \tensor h^{(2)}_\beta$ and suppressing the index
$\beta$).  As such we can form the semidirect product Hopf algebra
$U_{\A'}^\dagger(\lieh) \semidirect U_{\A'}^{\text{res}}(\lieg)$ (also
called the smash product, see Montgomery \cite{Montgomery93}).
This Hopf algebra is (densely) generated by
$\{E_i,F_i,K_i\} \cup \{\delta_\lambda\}_{\lambda \in \Lambda},$ with
the standard quantum group relations together with 
$$\delta_\lambda \delta_\gamma=\delta_{\lambda,\gamma}
\delta_\lambda,\qquad \sum_{\lambda \in \Lambda} \delta_\lambda=1$$
$$\delta_\lambda K_i=K_i \delta_\lambda \qquad \delta_\lambda
E_i= E_i \delta_{\lambda-\alpha_i} \qquad  \delta_\lambda F_i= F_i
\delta_{\lambda + \alpha_i}.$$

If $U_{\A'}^\dagger(\lieh) \semidirect U_{\A'}^{\text{res}}(\lieg)$
acts on an $\A'$-module $V,$ and $v\in V,$ we say $v$ is of weight $\lambda
\in \Lambda$ if $K_iv=q^{\bracket{\lambda,\alpha_i}}v$ and
$fv=f(\lambda)v$ for $f \in U_{\A'}^\dagger(\lieh),$ and we say $V$ is
a $\lambda$ weight space if it consists entirely of weight $\lambda$
vectors.  Let $\WW$ be the 
direct product of all $U_{\A'}^\dagger(\lieh) \semidirect
U_{\A'}^{\text{res}}(\lieg)$-modules which are a finite direct sum of
$\A'$-free 
$\lambda$ weight spaces for $\lambda \in \Lambda.$  Of course
$U_{\A'}^\dagger(\lieh) \semidirect U_{\A'}^{\text{res}}(\lieg)$  acts
on $\WW.$ The kernel of this action is a two-sided ideal  $I$ (clearly $I$
includes at least $K_i-q^{\alpha_i}$), and since
the tensor product of two finite direct sums of $\A'$-free
$\lambda$-spaces for $\lambda \in \Lambda$ is another such, $I$ is a
Hopf ideal.  Thus the quotient $U_{\A'}^\dagger(\lieh) \semidirect
U_{\A'}^{\text{res}}(\lieg)/I$ is a Hopf algebra which embeds into
$\End(\WW).$  The product topology on $\WW$ gives $\End(\WW)$
 a topology, one in which a sequence
converges if and only if it converges on each finite-dimensional
submodule.  The closure of 
$U_{\A'}^\dagger(\lieh) \semidirect 
U_{\A'}^{\text{res}}(\lieg)/I$ in this topology is what we call
$U_{\A'}^\dagger(\lieg).$  The product, coproduct, antipode and counit
clearly extend to the completion ($\Delta$ has range the closure of
$U_{\A'}^\dagger(\lieg) \tensor U_{\A'}^\dagger(\lieg)$ in
$\End_{\A'}(\WW \tensor \WW),$ which we will refer to as
$U_{\A'}^\dagger(\lieg) \,\overline{\tensor} \, U_{\A'}^\dagger(\lieg),$
the completed tensor product).

The reader might reasonably wonder, after the sequence of rather
abstract steps in the construction, whether there is anything
left to this algebra at all.  In fact $W^{\lambda}_{\A'}$ ($\lambda
\in \Lambda^+$) can be made
into a $U_{\A'}^\dagger(\lieh) \semidirect 
U_{\A'}^{\text{res}}(\lieg)$-module by letting $f\in
U_{\A'}^\dagger(\lieh)$ act on a weight $\lambda$ vector by
multiplication by $f(\lambda).$ $W^\lambda_{\A'}$ is  a finite direct sum of
free weight spaces, as above, so any pair of elements of the semidirect
product that act as different endomorphisms on some $W_{\A'}^\lambda$
represents different elements of $U_{\A'}^\dagger(\lieg).$

This extended Hopf algebra $U_{\A'}^\dagger(\lieg)$ is a ribbon Hopf
algebra (see \cite{CP94} for the definition). 
Specifically notice that our earlier $R$-matrix
$$q^{\sum_i\check{\alpha_i} \tensor
\lambda_i}$$
 is an element of
$U_{\A'}^\dagger(\lieg) \,
\overline{\tensor} \, U_{\A'}^\dagger(\lieg).$  Therefore so is
\begin{multline}\label{Rmatrix}R =
q^{\sum_i\check{\alpha_i} \tensor
\lambda_i} \sum_{t_1, \ldots t_N=1}^\infty
\prod_{r=1}^N q_{\beta_r}^{t_r(t_r+1)/2} (1-q_{\beta_r}^{-2})^{t_r}
[t_r]_{q_{\beta_r}}\!! E_{\beta_r}^{(t_r)} \tensor F_{\beta_r}^{(t_r)} \\ 
\in U^\dagger_{\A'}(\lieg) \,\overline{\tensor} \, U^\dagger_{\A'}(\lieg)
\end{multline}
where
$q_{\beta_r}=q^{d_i}$ when $\beta_r$ is the same length as
$\alpha_i.$   A fairly standard calculation confirms that $R$ is
a quasitriangular element for $U^\dagger_{\A'}(\lieg).$  Further,
notice that the grouplike element $q^\rho$ is a
charmed element of the Holf algebra for this $R$, making
$U^\dagger_{\A'}(\lieg)$ into a ribbon Hopf algebra.  In particular,
conjugation by $q^{2\rho}$ is the square of the
antipode, so that for any finite-dimensional $U^\dagger_{\A'}(\lieg)$
module $V,$ free over $\A',$ the functional 
$$\qtr_V \colon U^\dagger_{\A'}(\lieg) \to \A'$$
\[\label{qtr} \qtr_V(x)=
\tr_{V}(q^{2\rho }x)\] 
is  an invariant functional on
$U^\dagger_{\A'}(\lieg)$ in the sense that $\qtr(a^{(1)} b
S(a^{(2)}))=\epsilon(a) \qtr(b)$ (using Sweedler's notation as
above).

Define the \emph{quantum dimension}
$$\qdim(V)=\qtr_V(1)=\tr(q^{2\rho})$$
and in particular define $\qtr_\lambda=\qtr_{W^\lambda_{\A'}}$ and
$\qdim(\lambda)=\qdim(W^\lambda_{\A'}).$ 
Notice that since $\qtr_{V\tensor W}=\qtr_V\qtr_W,$ 
$$\qdim(V\tensor W)=\qdim(V)\qdim(W).$$
  Finally, a simple
calculation modeled on the classical Weyl character formula gives 
\begin{equation} \label{Weyl_character}
\qdim(\lambda)=\prod_{\beta >0} \left(q^{\bracket{\lambda+\rho,\beta}}
- 
q^{-\bracket{\lambda+\rho,\beta}}\right)/\left(q^{\bracket{\rho,\beta}}
- q^{-\bracket{\rho,\beta}}\right) \in \A. 
\end{equation}

\section{Roots of Unity and the Affine Weyl Group} \label{unity}

Now restrict the generic $q$ to a root of unity.  Specifically, let
$l$ be a positive integer, and consider the homomorphism $\A' \to
\QQ[\s],$ where $\s$ is an abstract primitive $lL$th root of unity
(i.e. satisfies the $lL$th cyclotomic polynomial) given by $s \mapsto
\s.$ As before define $U_\s^\dagger(\lieg)=U_{\A'}^\dagger(\lieg)
\tensor_{\A'} \QQ[\s],$ $W^\lambda_\s=W^\lambda_{\A'}\tensor_{\A'}
\QQ[\s].$ Write 
$\q=\s^L.$ Likewise $q \mapsto \q$ gives a homomorphism $\A \to
\QQ[\q],$ and $U_\A^{\text{res}}(\lieg) \to U_\q^{\text{res}}(\lieg).$
Write $\q^\lambda$ for $q^\lambda \tensor 1 \in U_{\A'}^\dagger(\lieg)
\tensor \QQ[\q]=U_\s^\dagger(\lieg).$

Notice $U_\s^\dagger \tensor U_\s^\dagger \iso (U_{\A'}^\dagger
\tensor U_{\A'}^\dagger) \tensor_{\A'} \QQ[\s]$ embeds naturally (and
densely in the inherited topology) into $(U_{\A'}^\dagger
\, \overline{\tensor} \, U_{\A'}^\dagger) \tensor_{\A'} \QQ[\s],$ and thus
we may define the latter space as the completed tensor product  $U_\s^\dagger
\,\overline{\tensor} \, U_\s^\dagger.$  $U_\s^\dagger$ then becomes a Hopf
algebra, and in fact a ribbon Hopf algebra since  the image of $R$ is in $U_\s^\dagger \,\overline{\tensor}\,
U_\s^\dagger.$

 For each $i \leq n$ let $l_i$ be $l/\gcd(l,d_i)$ (that is, the degree
 of $q_i$) and let $l_i'$ be $l_i$ or $l_i/2$ according to whether
 $l_i$ is odd or even (so that $l_i'$ is the least natural number such
 that $q_i^{l_i'} \in \{\pm 1\}.$ Likewise let $l'$ be $l$ or $l/2$
 according to whether $l$ is odd or even.   Define the \emph{affine Weyl group,}
 $\Weyl_l,$ to be the group of isometries of $\lieh^*$ generated by
 reflection about the hyperplanes
$$\bracket{x,\alpha_i}=\bracket{kl_i' \alpha_i/2,\alpha_i}=kl_i'd_i$$ for
each $k \in \ZZ$ and each $\alpha_i \in \Delta.$ This includes  the
Weyl group $\Weyl$ as a subgroup 
(when $k=0$).  Again we will usually be interested
in the translated action of the affine Weyl group, 
given by $\sigma \cdot \lambda=\sigma(\lambda+\rho) -\rho.$

\begin{lemma}\label{weyl_group}
The affine Weyl group is the semidirect product of the ordinary Weyl
with the group of
translations $l'\check{\Lambda}_r$ if $l'$ is  divisible by $D$
 or $l'\Lambda_r$ if $l'$ is not
divisible by $D.$ 
In particular a set of
generators consists of 
 reflections $\sigma_{\alpha_i},$ $\alpha_i
\in \Delta$ together with translation by $l'\theta/D$ (if $D|l'$), or
$l'\phi$ 
(if $D\!\! \not{|}l'$).  A fundamental domain for the translated
action of the affine Weyl
group is the \emph{principal Weyl alcove,} $C_{l},$ which is the region
$\bracket{x+\rho,\alpha_i} \geq 0,$ $\bracket{x+\rho,\theta} \leq l'$
(if $D|l'$) or $\bracket{x+\rho,\alpha_i} \geq 0,$
$\bracket{x+\rho,\phi} \leq l'$ otherwise.
\end{lemma}
\begin{pf}
Reflection about the hyperplane $\bracket{x,\alpha_i}=0$ followed by
reflection about 
$\bracket{x,\alpha_i}=\bracket{l_i'\alpha_i/2 ,\alpha_i}$ gives
translation by $l_i'\alpha_i,$ which since $d_i=1$ or $d_i=D$ and $D$ is prime
is translation by $l' \alpha_i$ or $l'\check{\alpha_i}$ according to the
divisibility of $l'.$   Conjugation by $\sigma_\beta$ for various
$\beta$ gives translation by $l'\gamma$ or $l'\check{\gamma}$ for
$\gamma$ any root.    Thus
$\Weyl_l$ contains the groups mentioned, and clearly is generated by
them.  Since the  Weyl group acts by conjugation on the
group of translations, the full group is a semidirect product. 

The subgroup of translations is generated by $l'\check{\beta}$
(resp.~$l'\beta$) for $\beta$ a long root of $\Phi$
(resp.~$\beta$ a short root of $\Phi$).  Thus a fundamental domain for
the group of translations would be the polygon bounded by the
hyperplanes $\bracket{x+\rho,\check{\beta}} \leq l'$
(resp.~$\bracket{x+\rho,\beta} \leq l'$) for all such $\beta.$ Since
this region is invariant under the translated action of $\Weyl,$ a
fundamental region for $\Weyl_l$ is given by the intersection of this
region with a fundamental region of this action of $\Weyl,$ which is
exactly the region given.
\end{pf}
 
\begin{remark} When $\lieg$ is not simply-laced the
affine Weyl group's action is distinctly different if $l'$
is divisible by $D$ versus if it is not.  When $l'$ is divisible by
$D,$ we recognize the action of the affine Weyl group described by
Kac \cite{Kac83} and many other authors discussing affine Lie algebras
and loop groups, except that the translations are
multiplied by $l'/D.$  This is the affine Weyl group relevant
to affine Lie algebras, and the affine Weyl group of the root
system $\Phi$ as discussed in, for example, Bourbaki
\cite{Bourbaki02}.  When $l'$ is not divisible by $D$ we recover
the affine Weyl group discussed in Jantzen \cite{Jantzen87} and many
other authors considering modular groups.  It is (with 
multiplication by $l'$) the usual affine Weyl group of the dual root system
$\check{\Phi}.$  It is remarkable that both versions
of the affine Weyl group appear, on equal footing, in the context of
quantum groups.
\end{remark}

\begin{lemma}\label{lm:different_Weyls}
The affine Weyl group above is the largest subgroup of $\Weyl^\dagger
\defequals \Weyl
\semidirect \frac{l}{2} \check{\Lambda}_r$ which fixes the root
lattice $\Lambda_r$ under the translated action.  In particular these
are equal when $2D | l.$
\end{lemma}

\begin{pf}
The subgroup of the group of translations $\frac{l}{2}
\check{\Lambda}_r$ which preserves $\Lambda_r$ is $\frac{l}{2}
\check{\Lambda}_r \cap \Lambda_r$=$l'\check{\Lambda}_r \cap
\Lambda_r$ since half a dual root is never a dual root.  If $D|l',$
then $ l' \check{\Lambda}_r \subset   
\Lambda_r.$   If $D$
does 
not divide $l'$ then they are relatively prime, so for short roots $l'
\check{\beta} =l' \beta \in \Lambda_r,$ but for long roots the smallest
multiple in $\Lambda_r$ is $D l' \check{\beta}=l'\beta,$ so
$\frac{l}{2} \check{\Lambda}_r \cap \Lambda_r \subset l'\Lambda_r.$
Since the translated action of the ordinary Weyl group preserves the
root lattice, the result follows.
\end{pf}
The remainder of the section is devoted to a proof of the quantum
version of Harish-Chandra's Theorem for $U_\s^\dagger(\lieg).$ 
Harish-Chandra's theorem for $U_{\A'}^\dagger(\lieg)$ (with of course
the classical Weyl group)  works by a
perfectly analogous argument (though simpler) and will
be left to the reader. 

Consider the action of the center $Z$ of $U_\s^\dagger$ on a Weyl
module of highest weight $\lambda.$  If $v$ is a vector of
weight $\lambda$ then so is $zv$ for any $z \in Z,$ so $z$ must act as
multiplication by an element of $\QQ[\s]$  on $v,$ say
$zv=\chi_\lambda(z)v.$ Since every element of the Weyl module
is of the form $Fv$ for some $F \in U_\s^{\dagger-}$ (the subalgebra of
$U_\s^\dagger $
generated by $\{1,F_i^{(k)}\}$) we have
$zFv=Fzv=\chi_\lambda(z) Fv,$ so that in fact $z$ acts as
multiplication by $\chi_\lambda(z)$ on the entire Weyl
module.  Thus each $\lambda$ gives us an algebra homomorphism
$\chi_\lambda$ from the
center $Z$ to $\QQ[\s].$

\begin{definition}If $\lambda,\gamma\in \Lambda$ say that $\lambda
\sim \gamma$ if $\chi_\lambda=\chi_\gamma.$
\end{definition}
Now $\chi_\lambda=\chi_\gamma$ if $\lambda$ occurs as a highest weight in
a highest weight $\gamma$ module, so $\sim$ includes at least the
extension of this inclusion relation to an equivalence relation. 

Let us first understand the relation $\sim$ in the case
$\lieg=\mathfrak{sl}_2.$ The Verma module of weight $j \theta,$ $j
\in \ZZ^{\geq 0}/2 \iso \Lambda^+,$ is a 
$U_\s^\dagger$-module spanned by $\{F^{(k)}v\,|\,k \geq 0\},$
where $v$ is of weight $j\theta$ and $E^{(k)}v=0$ for all $k>0.$
For this module
$$\left[\begin{matrix}2j+s-r\\s\end{matrix}\right]_\q E^{(r)}F^{(s)}v=
F^{(r-s)}  
v$$ 
if $s\leq r.$ 

Notice $[r]_\q=0$ if and only if $r$ is a
multiple of $l'.$ Thus
 $F^{(s)}v$ is a highest weight vector when
either $s=2j+1$ or $s=2j+1-kl'$ and $s<l'.$  So $j \theta  \sim
-(j+1) \theta$ when
$j \geq 0$ and $j\theta \sim (kl'-j-1) \theta$ when $2j<(k+1)l'.$  By transitivity $j
\sim j'$ whenever $j\theta$ is connected to $j'\theta$ by the affine
Weyl group.

Now consider a general $\lieg.$
\begin{proposition} \label{weyl_sim}
If $\lambda,\gamma \in \Lambda$ and there is a $\sigma \in \Weyl_l$
such that $\gamma=\sigma\cdot 
\lambda,$ then $\lambda \sim \gamma.$ 
\end{proposition}

\begin{pf}
Let $\lambda \in \Lambda.$ Recall that the Verma module of highest
weight $\lambda$ can be constructed as follows.  Consider
$U_\s^\dagger$ as a $U_\s^\dagger$-module
under the adjoint action, and quotient it by the left ideal generated
by $E_i^{(k)}$ and $K_i - q^{(\lambda, \alpha_i)})$ for all $i\leq N$
and $k \in \NN.$
It is easy to see that the vector $1$ (which above was called $v$) is
a highest weight vector of weight $\lambda.$ Now for each $i \leq N$
the set 
$\{E_i^{(k)},F_i^{(k)},K_i\}$ generate a subalgebra isomorphic to
$U_{s_i}^\dagger(\mathfrak{sl}_2)$ (here $\s_i=\s^{d_i}$)
and the vectors 
$F_i^{(k)}v$ span a $U_{\s_i}^\dagger(\mathfrak{sl}_2)$
module isomorphic to the Verma module of weight
$\bracket{\lambda,\check{\alpha_i}}/2.$ Therefore the vectors
$F_i^{\bracket{\lambda,\check{\alpha_i} }+1}v$ and
$F_i^{\bracket{\lambda,\check{\alpha_i} }+1-kl_i'}v$ where
$\bracket{\lambda,\check{\alpha_i} }<(k+1)l_i'$ give highest weight vectors.
Thus $\lambda \sim \sigma\cdot \lambda,$ for $\sigma$ a
generator of the affine Weyl group.  The result follows by the
transitivity of the $\sim$ relation.
\end{pf}

The other direction of Harish-Chandra's theorem requires the
$R$-matrix.

Writing $R=\sum_j u_j \tensor v_j$ (of course the sum is infinite),
define
$$D=\left(\sum_j v_j  \tensor u_j\right) \left(\sum_k u_k  \tensor
v_k\right)$$
in $U^\dagger_\s(\lieg) \,\overline{\tensor}\, U^\dagger_\s(\lieg).$  
In turn write $D=\sum_i x_i \tensor  y_i.$ Let
$$\Psi \colon (U^\dagger_\s(\lieg))^* \to U^\dagger_\s(\lieg)$$
$$\Psi(z^*)=\sum_i z^*(y_i) x_i,$$
where $(U^\dagger_\s(\lieg))^*$ is the direct sum over all $\lambda
\in \Lambda^+$ of the set of
functionals on $U^\dagger_\s(\lieg)$ which factor through the
representation on $W^\lambda_\s.$ 
This map is called the Drinfel'd map. We will also be interested in 
$$D_\lieh=\q^{\sum_i\lambda_i \tensor \check{\alpha_i}}
 \q^{\sum_i\check{\alpha_i} \tensor
\lambda_i}=\q^{2\sum_i\lambda_i \tensor \check{\alpha_i}}$$
and the associated 
$$\Psi_\lieh \colon (U^\dagger_\s(\lieg))^* \to U^\dagger_\s(\lieh)$$
$$\Psi_\lieh(z^*)=\sum_i z^*(y_i) x_i$$
writing $D_\lieh=\sum_i x_i \tensor y_i.$

By the PBW theorem  there is a
well-defined map
$$\Theta \colon U_\s^\dagger(\lieg) \to U_\s^\dagger(\lieh)$$
 given
by sending all products in the PBW basis which contain factors of
$E_i$ or $F_i$ to 
zero and all other products to themselves.  Thus $\chi_\lambda=\lambda
\circ \Theta$ on the center $Z.$  What's more, since the only terms in
$D$ which do not 
contain factors of the form $E$ and $F$ are those in $D_\lieh,$ 
\begin{equation} \label{theta_phi}
\Theta\Psi=\Psi_\lieh.
\end{equation}

Recall that the adjoint action of of $U_\s^\dagger$ on itself is
given by 
\[\ad_a(x)= a^{(1)} x S(a^{(2)})\]
where we have used Sweedler's notation which writes
$\Delta(a)=a^{(1)} \tensor a^{(2)}$ with understood summation sign
and indices.  An invariant element of $U_\s^\dagger$ is then an $x$
such that $\ad_a(x)=\epsilon(a) x$ for all $a\in U_\s^\dagger.$ We
note 
\cite{Kuperberg91} argues that for any 
Hopf algebra the map $(a,b) \mapsto (a^{(1)},ba^{(2)})$ is one-to-one
and onto from $U_\s^\dagger \,\overline{\tensor} \, U_\s^\dagger$ to itself
(his argument was for ordinary 
tensor product, but is easily adapted to our topological situation).
For a given $u \in U_\s^\dagger$ choosing $(a,b)$ which is mapped to
$(u,1),$ one 
argues that for an ad-invariant element $z,$ we have  $uz$=$a^{(1)} z S(a^{(2)})
S(b)= \epsilon(a) z S(b)= z a^{(1)} S(a^{(2)}) S(b)=zu$ for all $u.$
Since the converse is clear we conclude that the invariant elements
of $U_\s^\dagger$ (indeed of any Hopf algebra) are exactly the
elements of the center.

Likewise the coadjoint action on $(U_\s^\dagger)^*$ sends $z^*$ to
\[\coad_a(z^*)=z^*(a^{(1)} \,\cdot\, S(a^{(2)})).\]

\begin{lemma}
The Drinfel'd map $\Psi$ takes invariant functionals to the center of
$U_\s^\dagger.$ 
\end{lemma}

\begin{pf}
  For notational
convenience  write
$$\Delta^2(a)= b^{(1)} \tensor  b^{(2)} \tensor b^{(3)}$$
and 
$$\Delta^3(a)= c^{(1)} \tensor c^{(2)} \tensor c^{(3)}
\tensor c^{(4)}.$$

Notice by the quasitriangularity of $R$
$$\Delta(a) D=D\Delta(a)$$
for all $a \in U_\s^\dagger.$  Then if $z^*$ is an invariant functional
\[\epsilon(a) \Psi(z^*)=  \sum_i z^*(y_ic^{(2)}S(c^{(3)})) x_i c^{(1)}
S(c^{(4)})\]
by the basic relations of a Hopf algebra so 
\begin{eqnarray*}
&= &\sum_i z^*(c^{(2)} y_i S(c^{(3)}))
c^{(1)} x_i S(c^{(4)})\\
=\sum_i \ad_{S^{-1}(b^{(2)})}(z^*)(y_i) b^{(1)} x_i S(b^{(3)})\\
=\sum_i \epsilon(S^{-1}(b^{(2)})) z^*(y_i) b^{(1)} x_i S(b^{(3)})\\
=\sum_i z^*(y_i) a^{(1)} x_i S(a^{(2)}).
\end{eqnarray*}
Thus $\Psi(z^*)$ is an invariant element of $U^\dagger_\s.$

\end{pf}

 \begin{corollary} If $\lambda \sim \gamma,$ then $\lambda  \Psi_\lieh=
 \gamma\Psi_\lieh$ on invariant functionals.
\end{corollary}

 \begin{proposition} \label{pr:almost_hc} Suppose $\lambda, \gamma \in
 \frac{1}{2L} \check{\Lambda}_r.$  If $\lambda \Psi_{\lieh}$
 agrees with $\gamma \Psi_{\lieh}$ on   	
 invariant functionals, then $\lambda$ and $\gamma$ are in the same
 orbit of $\Weyl^\dagger.$    
  Further, the set $\{\lambda\Psi_\lieh\},$ where
 $\lambda$ runs through a choice of representative of each equivalence
 class in $\frac{1}{2L} \check{\Lambda}_r/\Weyl^\dagger,$
 is a set of linearly independent
 functionals on the quantum traces.  
\end{proposition}

\begin{pf}
For each $\nu \in \Lambda^+ $ the functional $\qtr_\nu$ is an
invariant functional.  By induction on the ordering 
we can form a linear 
combination of these $\qtr_\nu$ to produce an invariant functional
which on $U_\s^\dagger(\lieh)$ acts as $\sum_{\sigma \in \Weyl}
\sigma(\nu)(\q^{2\rho} \,\cdot\,)$ for each $\nu \in \Lambda^+.$  Notice
that $(\lambda \tensor 
\mu)(D_\lieh)=\q^{2\bracket{\lambda,\mu}}$  for $\mu \in \Lambda.$  Thus
\[
\lambda(\Psi_\lieh(\sum_{\sigma \in \Weyl}\sigma(\nu)(\q^{2\rho}
\,\cdot\,)))
=(\lambda \tensor \sum_{\sigma \in \Weyl}\sigma(\nu)\left(D_\lieh (1
\tensor \q^{2\rho})\right))
=\sum_{\sigma \in \Weyl} \q^{2\bracket{\lambda+\rho,\sigma(\nu)}}.
\]
The set of maps $\{\q^{2\bracket{\lambda+\rho,\,\cdot\,}}: \lambda \in
\frac{1}{2L}\check{\Lambda}_r/\frac{l}{2}
\check{\Lambda}_ r\}$ is a
basis for maps from $\Lambda/(lL\Lambda)$ to $\QQ[\s].$  $\Weyl$
permutes this basis, so $\sum_{\sigma \in \Weyl}
\q^{2\bracket{\lambda+\rho,\sigma(\,\cdot\,)}}$ forms a basis for maps
from $ \left(\Lambda/(lL\Lambda)\right)^{\Weyl}$ to $\QQ[\s],$ when
$\lambda$ ranges over 
representatives of each Weyl orbit in
$\Weyl^\dagger.$  Thus $\lambda\Psi_\lieh$ is unchanged by $\Weyl
\semidirect \frac{l}{2} \check{\Lambda}_r,$ and any set of orbit
representatives is linearly independent.
\end{pf}

\begin{lemma} \label{root_center} If $\lambda \sim \gamma$ then
$\lambda-\gamma \in \Lambda_r.$
\end{lemma}

\begin{pf}
Notice an element $f \in U_\s^\dagger(\lieh)$ is in the center of
$U_\s^\dagger(\lieg)$  if and only if $f(\lambda)=f(\lambda+\alpha_i)$
for all $\lambda \in \Lambda$ and all $i.$  The sub-Hopf algebra of
such functions is isomorphic to the Hopf algebra of functions on the
fundamental group $\Lambda/\Lambda_r.$  Such an $f$ acts on $\lambda$
by multiplication by $f(\lambda),$ so every such $f$ will agree on
$\lambda$ and $\gamma$ if and only if $\lambda-\gamma \in \Lambda_r.$

\end{pf}

 \begin{theorem} \label{harish_chandra} $\lambda \sim \gamma$ if and
 only if $\lambda=\sigma 
 \cdot \gamma$ for some $\sigma \in \Weyl_l.$
\end{theorem}

\begin{pf}
That the latter implies the former is exactly Proposition
\ref{weyl_sim}.

If $\lambda \sim \gamma,$ then by  Proposition \ref{pr:almost_hc} they
are connected by an element of $\Weyl^\dagger.$  On the other hand by
Lemma \ref{root_center} they differ by an element of the root
lattice, and thus the element of $\Weyl^\dagger$ must preserve the
root lattice (it is easy to see that if an element of $\Weyl^\dagger$
takes one vector to another vector that differs from it by a root
vector, the difference of any vector and its image is a root
vector).  Thus by Lemma \ref{lm:different_Weyls} they are connected
by an element of the affine Weyl group $\Weyl_l.$

\end{pf}

\begin{corollary}\label{independent}
In fact, $\{\chi_\lambda\},$
choosing one $\lambda$ from 
each translated $\Weyl_l$ equivalence class, is linearly independent
as a set of functionals on the center.
\end{corollary}

\begin{pf} By Proposition \ref{pr:almost_hc} a linear relation between
these would reduce to a linear relation between those in the translated
$\Weyl^\dagger$ orbit of some
$\lambda.$ Since elements of this orbit which are not $\Weyl_l$
equivalent must be in distinct classes of $\Lambda/\Lambda_r,$
computing them on the center intersected with $U_\s^\dagger(\lieh)$
shows that no such nontrivial relation exists.
\end{pf}

\section{Weyl Filtrations and Tilting Modules}

In this section use $U$ to refer to any of the forms of the quantum
group defined in the last two sections: $U_q(\lieg),$
$U_\A^{\text{res}}(\lieg),$ 
 $U_{\A'}^\dagger(\lieg),$  $U_{\q}^{\text{res}}(\lieg)$ or
 $U_{\s}^\dagger(\lieg),$ and use ``the ground ring'' to refer to
 $\QQ(q),$ $\A,$ $\A',$ $\QQ[\q],$ or $\QQ[\s]$ as appropriate.  Also,
 drop the
subscripts from such notation as $W_{\A'}^\lambda$ when no confusion
would ensue.  While many of the results of this section will apply to
all forms of the quantum group, we will be interested in their
application only to $U_\q^{\text{res}}$ and $U_\s^\dagger.$  

\begin{definition} a $U$-module $V$  is said to have a \emph{Weyl filtration}
if there exists a sequence of submodules 
$$\{0\}=V_0 \subset V_1 \subset \cdots \subset V_{n-1} \subset V_n=V$$
such that for each $1 \leq i \leq n,$ $V_i/V_{i-1}$ is isomorphic to
the Weyl module $W^{\lambda}$ for some $\lambda \in \Lambda^+.$
\end{definition}

\begin{proposition} \label{WF}
Suppose $W$ is a $U_{\A}^{\text{res}}(\lieg)$-module such that $W
\tensor_{\A} \QQ(q) =\bigoplus_i W^{\lambda_i}_q.$  Then $W \tensor_\A
\QQ[\q]$ and $(W \tensor_\A \A') \tensor_{\A'} \QQ[\s]$ admit Weyl
filtrations with the $i$th factor of highest weight $\lambda_i,$ where
the $\lambda_i$ are assumed to be ordered so that $\lambda_i$ is never
greater than $\lambda_j$ for $j<i.$ 
\end{proposition}

\begin{pf}  We shall prove the proposition over
$\QQ[\q],$ the argument is exactly the same for
$U_\s^\dagger(\lieg).$ 

Decomposing $W=W_{\text{tor}} \oplus W_{\text{free}}$ into its torsion
and free parts over $\A,$ notice that $W \tensor_{\A} \QQ(q)=W_{\text{free}}
\tensor_{\A} \QQ(q)$ and likewise for $\QQ[\q].$  Since
$W_{\text{tor}}$ is a $U_{\A}^{\text{res}}$-module, the
quotient by it is a free $\A$-module and a
$U_{\A}^{\text{res}}$ module whose tensor with $\QQ(q)$ and
$\QQ[\q]$ are isomorphic to that of $W.$  Thus we can assume $W$ is a free
$\A$-module.  Notice in this case the maps $v \mapsto v \tensor 1$ are
 injective maps from $W$ to $W \tensor_\A \QQ(q)$ and $W \tensor_\A
 \QQ[\q]$ whose ranges span.

Let $w \in W$ be such that $w \tensor 1 \in W \tensor_{\A} \QQ(q)$ is
a vector of weight $\lambda_1.$ By the maximality of $\lambda_1$ $w$
must be a highest weight vector.  Then
$U_\A^{\text{res}}w$ is a $U_\A^{\text{res}}$-module,
free over $\A,$ whose tensor product with $\QQ(q)$ yields a
$U_q$-module isomorphic to $W^{\lambda_1}_q.$  Thus
$U_\A^{\text{res}}w$ must be 
isomorphic to 
$W^{\lambda_1}_\A.$  Its tensor product with $\QQ[\q]$  gives a
submodule isomorphic to $ W^{\lambda_1}_\q.$  Therefore the quotient
$W/W^{\lambda_1}_\A$ is a module whose tensor product with $\QQ(q)$ is
isomorphic to $\bigoplus_{i>1} W^{\lambda_i}_q$ and whose tensor
product with $\QQ[\q]$ is $(W \tensor_\A \QQ[\q])/W^{\lambda_1}_\q.$
By induction the proposition follows.
\end{pf}

\begin{corollary}\label{tensor_product}
The tensor product of two $U_\q^{\text{res}}$  or $U_\s^\dagger$
modules with a Weyl filtration admits a Weyl 
filtration.
\end{corollary}

\begin{pf} By induction it suffices to
prove that $W^\lambda \tensor W^\gamma$ admits a Weyl filtration.  This
follows from the previous proposition.  
\end{pf}
\begin{remark}  Notice the entries in that Weyl
filtration are the same as the entries in the classical decomposition of
the tensor product of classical modules which were direct sums with
the same entries as the original Weyl filtrations.  Thus if we
restrict attention to  modules with a Weyl filtration the category of
such modules forms a monoidal category which is ``the same'' as the
tensor category of classical finite-dimensional modules if we replace
the notion of direct sum decomposition of modules with that of Weyl
filtration. 
\end{remark}

Recall if $V$ is a module over $U$ the space of linear functionals
from $V$ to the ground ring is naturally a right $U$-module $V^*.$ We
can compose the induced representation with either $S$ or $S^{-1}$ to make it
a left module which we also call $V^*$ (the two ways of doing this
give distinct but isomorphic module structures, and the distinction
will not be relevant to us). If $V$ is a finite direct sum of free
finite-rank weight spaces, so is $V^*.$ In this case, $V$
is isomorphic as a $U$-module to $V^{**}$ (though not by the canonical
identification 
between these two as modules over the ground ring, instead this
identification must be 
composed with conjugation by $q^{2\rho}$).  The dual of the tensor
product of two modules is isomorphic to the tensor product of the
duals.

\begin{definition} a $U$-module $V$ is said to have a \emph{dual Weyl
filtration} 
if there exists a sequence of submodules 
$$\{0\}=V_0 \subset V_1 \subset \cdots \subset V_{n-1} \subset V_n=V$$
such that for each $1 \leq i \leq n $ $V_i/V_{i-1}$ is isomorphic to
the dual of a Weyl module $(W^{\lambda})^*$ for some $\lambda \in
\Lambda^+.$
\end{definition}

Note that $V$ admits a dual Weyl filtration if and only if $V^*$ admits a
Weyl filtration.

\begin{definition}  A $U$-module $V$ is a \emph{tilting module} if it
admits both a Weyl filtration and a dual Weyl filtration. 
\end{definition}

\begin{corollary}\label{tilting} The properties of admitting a Weyl filtration,
admitting a good filtration or being a tilting module are preserved by
tensor product.   
\end{corollary}

\begin{remark}
The category of tilting modules forms a ribbon category
(i.e. with tensor products and duals) which is not semisimple, but
because of the existence of Weyl and dual Weyl filtrations, behaves in
many respects like the semisimple tensor category of classical $\lieg$
modules, in particular as far as the link invariant is concerned.
\end{remark}

If $V\iso W \bigoplus W',$ then $V$ is tilting if and only if $W$ and
$W'$ are tilting, so to understand tilting modules it suffices to
understand indecomposable tilting modules.  To do this requires a
short detour into elementary homological algebra.  For more detail on
the subject, see Mac Lane \cite{MacLane63}.  

Recall that for two modules $A$ and $C$ over a ring, the set of exact
sequences 
$$0 \to C \to ? \to \cdots \to ? \to A \to 0$$
forms a chain complex indexed by the number of intervening modules,
with an appropriate boundary operator.  The associated homology
$\Ext^n(A,C)$ is functorial in each variable, and if 
$$
0 \to X \to Y \to Z \to 0$$
is a short exact sequence then we get the long exact sequences of homology
\begin{equation} \label{lesh1} \begin{CD}
0 @>>> \Ext^0(Z,C) @>>> \Ext^0(Y,C) @>>> \Ext^0(X,C) @>>> \\
& & \Ext^1(Z,C) @>>> \Ext^1(Y,C) @>>> \cdots
\end{CD}
\end{equation}
and
\begin{equation} \label{lesh2} \begin{CD}
0 @>>> \Ext^0(A,X) @>>> \Ext^0(A,Y) @>>> \Ext^0(A,Z) @>>> \\
& & \Ext^1(A,X) @>>> \Ext^1(A,Y) @>>> \cdots.
\end{CD}\end{equation}
Finally, $\Ext^0(A,C) \iso \Hom(A,C)$ and if $A$ and $C$ are such that
every short exact sequence $0 \to C \to B \to A \to 0$ splits, then
$\Ext^1(A,C)=0.$  

\begin{lemma} \label{splits}
If $V$ is a $U$-module and $\lambda \in \Lambda^+$ is such that no
weight appearing in $V$ is greater than $\lambda$ then any quotient
$$\begin{CD} V @>g>> W^\lambda @>>> 0
\end{CD}$$
splits.
\end{lemma}

\begin{pf}
There must be a vector $v \in V$ of weight $\lambda$ which is in the
preimage of a highest weight vector in $W^\lambda$ under $g.$ By the
condition on $\lambda,$ $v$ is of maximal weight in $V$ and hence is a
highest weight vector, so there is a map $g':W^\lambda \to V$ sending
a highest weight vector to $v.$ Clearly $gg'$ is nonzero on the
highest weight vector and hence is a multiple of $1,$ so the sequence
splits.
\end{pf}

\begin{lemma} \label{ext}
If $A$ admits a Weyl filtration, and $C$ admits a dual Weyl filtration then
$\Ext^1(A,C)=0.$
\end{lemma}

\begin{pf}
We will prove first the base case, then do induction on the filtration
of $A,$ then induction on the filtration of $C.$
\begin{itemize}
\item \emph{$\Ext^1(W^\lambda, (W^\gamma)^*)=0$ for all $\lambda,\gamma \in
      \Lambda^+.$} Suppose first that $\lambda \not \! < \gamma^*,$
      where $\gamma^*,$ which is minus the image of $\gamma$ under the
      action of the 
      longest element of the Weyl group, is the maximal highest weight
      of $(W^\gamma)^*.$ Then if
$$ 0 \to (W^\gamma)^* \to B \to W^\lambda \to 0
$$
 the sequence splits by Lemma
\ref{splits} and the result follows. 
  On the other hand if $\lambda <\gamma^*$ and
$$0 \to (W^\gamma)^* \to B \to W^\lambda \to 0$$ then dualizing
$$0 \to (W^\lambda)^* \to B^* \to W^\gamma \to 0$$
and again by Lemma \ref{splits} the result follows.
\item \emph{$\Ext(A,(W^\gamma)^*)=0$ if $A$ admits a Weyl filtration.}
      By induction there is a short exact sequence
$$0 \to W^\lambda \to A \to A' \to 0$$ with $A'$ admitting a Weyl
filtration and hence $\Ext^1(A',(W^\gamma)^*)=0.$ From the long exact
sequence (\ref{lesh1})
$$\to [\Ext^1(W^\lambda,(W^\gamma)^*)=0] \to \Ext^1(A,(W^\gamma)^*) 
\to [0=\Ext^1(A',(W^\gamma)^*)] \to $$
from which it follows $\Ext(A,(W^\gamma)^*)=0.$
\item \emph{$Ext(A,C)=0$ if $A$ admits a Weyl filtration and $C$ admits a
      dual Weyl filtration.}  Again inductively we have a sequence
$$0 \to (W^\gamma)^* \to C \to C' \to 0$$ with $C'$ admitting a dual Weyl
filtration and hence $\Ext^1(A,C)=0.$ Again the long exact sequence
\ref{lesh2} and the previous two items give $\Ext(A,C)=0.$
\end{itemize}
\end{pf}

\begin{proposition}\label{unique}
If $Q,$ $Q'$ are indecomposable tilting modules over $U$ each with a maximal
vector of weight $\lambda,$ then $Q \iso Q'.$
\end{proposition}

\begin{pf}
Suppose $v$ is a weight $\lambda$ vector in $Q$ and $v'$ is a weight
$\lambda$ vector in $Q'.$  Let  $f:W^\lambda \to Q$ and $f':W^\lambda
\to Q'$ send a particular highest weight vector to $v$ and $v'$
respectively.  Let $j$ be the smallest integer such that $V_j$ contains $v$ in
a Weyl filtration of $Q.$ Then $V_j/V_{j-1} \iso W^\lambda.$ By
the maximality of $\lambda,$ $V_{j-1}/V_{j-2} \iso W^\gamma$ with
$\lambda \not\!< \gamma.$ By Lemma \ref{splits} we can find a new
$V_{j-1}'$ (without changing $V_{j-2}$) such that the filtration is
still Weyl but $v$ is now an element of $V_{j-1}.$ Inductively, there
exists a Weyl filtration with the image of $f$ being $V_1,$ which is
to say there is a short exact sequence
$$\begin{CD}
0 @>>> W^\lambda @>f>> Q @>>> N @>>> 0
\end{CD}$$
with $N$ admitting a Weyl filtration.   The long exact sequence
(\ref{lesh1}) gives  
$$0 \to \Hom(N,Q') \to \Hom(Q,Q') \to \Hom(W^\lambda,Q')
 \to [0=\Ext^1(N,Q') ]\to$$
by Lemma \ref{ext}, so that the sequence is in fact short exact, and
$f' \in \Hom(W^\lambda, Q')$ must factor through a map $g': Q \to Q'$
which takes $v$ to a 
nonzero multiple of $v'.$ By the same argument with $Q$ and $Q'$
reversed there is a map $g:Q' \to Q$ taking $v'$ to a nonzero multiple
of $v.$  Thus $gg'$ is a map from $Q$ to itself taking $v$ to a
nonzero multiple of itself.

$\{(gg')^n[Q]\}_{n \in \NN}$ is a nested sequence of submodules and
thus by finite-dimensionality must stabilize on some submodule
$(gg')^M[Q]$ such that $gg'$ is onto when restricted to this submodule.
So $Q=(gg')^M[Q] \bigoplus \Ker((gg')^M).$  Since $Q$ is
indecomposable one of these summands must be zero, and since $v \in
(gg')^MQ$ one has $Q=(gg')^M[Q]$ so $gg'$ is invertible. 
\end{pf}
\begin{corollary}
Eery tilting module is a direct sum of indecomposable tilting
modules.  Every indecomposable tilting module is isomorphic to some
$T_\lambda,$ the unique indecomposable indecomposable tilting module
with a maximal vector of  weight $\lambda.$
\end{corollary}

\begin{corollary}[Linkage Principle] \label{linkage} A simple module with
highest weight $\lambda$ can occur as a composition factor in the
Weyl or indecomposable tilting module of highest weight $\gamma$ only
if $\lambda \leq \gamma$ 
and $\lambda=\sigma\cdot \gamma$ for some $\sigma \in \Weyl_l.$ 
\end{corollary}

\section{Negligible Modules and the Weyl Alcove}

As in Section \ref{unity} we work with $\s$ a primitive $lL$th root of
unity, and consider the quantum group $U_\s^\dagger(\lieg).$ Let
$M$ be the lattice of translations $ l'\check{\Lambda}_r$ or $ l' \Lambda_r$ according to
whether $D$ divides $l'$ or not, so that $\Weyl_l=\Weyl \semidirect
M.$ Define hyperplanes
$$w_{k,\alpha}=\begin{cases}
\{x \in \lieh^*\,,\,\bracket{x + \rho,\alpha}=k l'\} & \text{if $D|l$}\\
\{x \in \lieh^*\,,\,\bracket{x + \rho,\check{\alpha}}=k l'\} &
\text{else,}
\end{cases}
$$
for $\alpha \in \Phi^+,$ called the \emph{walls} of the (translated)
affine Weyl 
group, so that the translated action of $\Weyl_l$ is generated by
reflections $\sigma_{k,\alpha}$ about $w_{k,\alpha}.$ 
These hyperplanes 
divide $\lieh^*$ into compact regions, called \emph{alcoves,} including the
principal alcove $C_l,$ such that $\sigma \mapsto \sigma\cdot C_l$
is a bijection between elements of $\Weyl_l$ and alcoves.  The
hyperplanes are called the \emph{walls} of the alcoves. Likewise the
\emph{walls} of $\Weyl \semidirect \frac{l}{2}\check{\Lambda}_r$ are
the hyperplanes $\bracket{x+\rho,\alpha}=kl/2.$ If $x$ is on the wall
$w_{k,\alpha}$ and no other wall then the stabilizer of $x$ in
$\Weyl_l$ is $\{1,\sigma_{k,\alpha}\}.$ A wall $w_{k,\alpha}$ of
an alcove $\sigma\cdot C_l$ is called a \emph{lower wall} if every point
$y$ in the interior satisfies $\bracket{y+\rho,\alpha}$ is greater
than the corresponding quantity for the points on the wall, and an
\emph{upper wall} otherwise.  Finally, let $\theta_0$ be $\theta$ if $D$
divides $l'$ and $\phi$ otherwise, so
that $w_{1,\theta_0}$ is the unique upper wall of $C_l,$ and the
intersection of the interior of $C_l$ with $\Lambda$ is 
$$\Lambda^l \defequals \{\lambda \in \Lambda^+ \, |\,
\bracket{\lambda+\rho,\theta_0}<l'\}.$$ 

A module $V$ is called \emph{negligible} if every intertwiner $\phi
\colon  \to V$ has quantumm trace $0.$

\begin{theorem} \label{negligible} In $U_\s^\dagger(\lieg),$ every
$T_\lambda$ with 
$\lambda$ not in $\Lambda^l$ is negligible, provided $l'\geq
D\check{h}$ if $D|l'$ or $l'> h$ otherwise.
\end{theorem}

\begin{lemma}
an indecomposible titling module is negligible if and only if its
quantum dimension is zero.
\end{lemma}

\begin{pf}
The algebra of intertwiners from an indecomposable module to itself
consists of multiples of the identity and nipotent intertwiners.
Since intertwiners commute with $q^{2\rho},$ they hae quantum trace
zero.  If the module has quantum dimension zero, then multiples of
the identity hae quantum trace zero.
\end{pf}

\begin{lemma}
The tensor product of a neglible module and another module is
a direct sum of negligible indecomposable modules.
\end{lemma}

\begin{pf}
If $V$ is a negligible module and $W$ is not, then the intertwiners
on $V \tensor W$ are just the tensor product of the set of
intertiners on $V$ tensored with those on $W.$  Thus the quantum
trace of any such intertwiner is the sum of a product of quantum
traces of intertwiners on $V$ and $W.$ In each term the first entry
of the product is zero.

If a module is negligible, so are all its direct summands.
\end{pf}

\begin{lemma}\label{walls} \hspace{1in}
\begin{alist}
\item If $\sigma \in \Weyl^\dagger,$ then
      $\qdim(\lambda)=(-1)^\sigma \qdim(\sigma\cdot \lambda)$ whenever
      $\lambda, \sigma  \cdot \lambda \in \Lambda^+,$ where $(-1)^\sigma$
      represents the orientation of $\sigma.$  In particular this is true of
      $\sigma \in \Weyl_l.$  Further, $\qdim(\lambda)=0$ if and only
      if $\lambda$ has  nontrivial stabilizer in $\Weyl \semidirect \frac{l}{2} \check{\Lambda}_r.$  
\item Every $T_\lambda$ where $\lambda$ is on a wall of $\Weyl_l$
  is negligible.
\end{alist}
\end{lemma}
\begin{pf} \hspace{1in}
\begin{alist}
\item 
By the Weyl formula
(\ref{Weyl_character}), $\qdim(\lambda)$ is 
$$
\qdim(\lambda)=\prod_{\beta >0} \left(\q^{\bracket{\lambda+\rho,\beta}}
- 
\q^{-\bracket{\lambda+\rho,\beta}}\right)/\left(\q^{\bracket{\rho,\beta}}
- \q^{-\bracket{\rho,\beta}}\right).$$
In fact we can interpret $\qdim(\lambda)$ by this formula even when
$\lambda$ is not in $\Lambda^+.$  It suffices to prove the first sentence
when $\sigma$ is a generator of the classical Weyl group
$\sigma_{0,\alpha_i}=\sigma_{\alpha_i}$ and when $\sigma$ is translation by
$l\check{\theta}/2.$  

Suppose first that $\sigma$ is $\sigma_{\alpha_i},$ then
\begin{multline*}\qdim(\sigma\cdot\lambda) = \prod_{\beta >0}
\left(\q^{(\sigma\cdot\lambda + \rho,\beta)} - \q^{-(\sigma\cdot \lambda +
\rho,\beta)}\right) / \left(\q^{(\rho, \beta)} -
\q^{-(\rho,\beta)}\right)  
\\= \prod_{\beta >0} \left(\q^{(\sigma_{\alpha_i}(\lambda +
\rho),\beta)} - \q^{(-\sigma_{\alpha_i}(\lambda + \rho),\beta)}\right)
/ \left(\q^{(\rho, \beta)} - \q^{-(\rho,\beta)}\right)  
\\= \prod_{\beta >0} \left(\q^{(\lambda +
\rho,\sigma_{\alpha_i}(\beta))} - \q^{-(\lambda +
\rho,\sigma_{\alpha_i}(\beta))}\right) / \left(\q^{(\rho, \beta)} -
\q^{-(\rho,\beta)}\right) 
\end{multline*}
since $\sigma_{\alpha_i}$ is a unipotent isometry.  Notice that
$\sigma_{\alpha_i}$ permutes the positive roots of $\Phi$ except for
$\alpha_i$ which it reverses (\cite{Humphreys72}[10.2]) so all factors
above stay the same except for one which changes sign.  Thus the
formula above gives $-\qdim(\sigma \cdot
\lambda)=(-1)^{\sigma}\qdim(\lambda).$  

Now suppose $\sigma$ is translation by $l\check{\theta}/2.$ Then
\begin{multline*}\qdim(\sigma\cdot \lambda) = \prod_{\beta >0}
\left(\q^{(\sigma\cdot \lambda + \rho,\beta)} - \q^{-(\sigma\cdot \lambda +
\rho,\beta)}\right) / \left(\q^{(\rho, \beta)} -
\q^{-(\rho,\beta)}\right)  
\\= \prod_{\beta >0} \left(\q^{(l\check{\theta}/2+ \lambda + \rho,\beta)} -
\q^{-(l\check{\theta}/2 +\lambda + \rho,\beta)}\right) /
\left(\q^{(\rho, \beta)} - \q^{-(\rho,\beta)}\right)  
\\= \prod_{\beta >0} \left(\q^{l\bracket{\beta,\check{\theta}}/2}\q^{(\lambda +
\rho,\beta)} - \q^{-l\bracket{\beta,\check{\theta}}/2}\q^{-(\lambda +
\rho,\beta)}\right) / \left(\q^{(\rho, \beta)} -
\q^{-(\rho,\beta)}\right) 
\\=\left(\prod_{\beta >0}
\q^{l\bracket{\beta,\check{\theta}}/2}\right)\prod_{\beta >0}
\left(\q^{(\lambda + \rho,\beta)} - \q^{-(\lambda 
+ \rho,\beta)}\right) / \left(\q^{(\rho, \beta)} -
\q^{-(\rho,\beta)}\right)\\  
=\q^{l\bracket{\rho,\check{\theta}}}\qdim(\lambda)\\
=\qdim(\lambda). 
\end{multline*}

Since the affine Weyl group is a subgroup of $\Weyl^\dagger,$ the same
result applies to the smaller 
group.

Of course if  $\lambda$ has nontrivial stabilizer than it lies on some
wall   so there is a reflection $\sigma$ which fixes $\lambda$ and  
      $\qdim(\lambda)=\qdim(\sigma \cdot
      \lambda)=-\qdim(\lambda),$ so $\qdim(\lambda)=0.$  $\lambda$ has
      no stabilizer
      $\qdim(\lambda)$ is a product of nonzero quantities, and thus nonzero.  
\item By the Linkage Principle, Corollary \ref{linkage}, $T_\lambda$
      has a Weyl filtration all of whose entries are affine Weyl
      equivalent to $\lambda.$  If $\lambda$ is on a wall, so are all
      weights in its affine orbit, and hence the quantum dimension of
      $T_\lambda,$ which is a sum of the quantum dimensions of the
      entries of the Weyl filtration, is zero.  Therefore $T_\lambda$
      is negligible.
\end{alist}  
\end{pf}

\begin{corollary} \label{qtr_cor} As functionals on the center
$$\qtr_{\sigma\cdot \lambda}=(-1)^\sigma \qtr_\lambda$$
when $\sigma \in \Weyl_l.$  As functionals on the image of quantum
traces under $\Psi$ the same is true when $\sigma \in \Weyl
\semidirect \frac{l}{2} \check{\Lambda}_r$  and $\sigma\cdot \lambda
\in \Lambda.$
\end{corollary}

\begin{lemma}\label{tensor_factors} \hspace{1in}\begin{alist}

\item If $W_\mu$ appears in a Weyl filtration of a tilting module $T,$ and no
      $\mu'$ in the translated $\Weyl_l$ orbit of $\mu$ with
      $\mu'>\mu$ appears 
      in that filtration, then 
      $T_\mu$ is a direct summand of $T.$ 
\item If $\mu$ appears as a highest weight in the \emph{classical} module
      $W^\lambda \tensor W^\gamma$ and no $\mu'>\mu$ in the translated
      $\Weyl_l$
      orbit of $\mu$ appears in the classical module $W^{\lambda'}
      \tensor W^{\gamma'}$ for $\lambda'\leq \lambda$ and $\gamma'\leq
      \gamma$
      in the translated $\Weyl_l$ orbits of $\lambda$ and $\gamma$ respectively,
      then $T_\mu$ is a direct summand of $T_\lambda \tensor T_\gamma.$
\end{alist}
\end{lemma}
\begin{pf} \hspace{1in}\begin{alist} 
\item By the Linkage Principle, Corollary \ref{linkage}, $W_\mu$
      appears in the filtration of an indecomposable direct summand 
      whose Weyl decomposition contains only modules with highest
      weights in the translated 
      $\Weyl_l$ orbit of $\mu.$  By the assumption on $\mu$ the weight
      $\mu$ is maximal in this summand, which must thus be
      isomorphic to $T_\mu$ by Proposition \ref{unique}. 
\item  Of course a factor of  $W_\mu$  must appear in a filtration of 
       $T_\lambda \tensor T_\gamma$ by Proposition~\ref{WF} and
       Corollary~\ref{tensor_product}. If a larger $\mu'$ in the orbit
       of $\mu$ also appeared in the filtration, it would appear in
       the classical decomposition of some $W^{\lambda'}
      \tensor W^{\gamma'}$ with $\lambda'$ and $\gamma'$ in a Weyl
      filtration of $T_\lambda$ and $T_\gamma$ respectively.  This is
      ruled out by the assumption, so by part (a) we are done.
\end{alist} 
\end{pf}

\begin{lemma}\label{wall_factor}  Suppose $\lambda, \gamma, \lambda+\sigma(\gamma) \in
\Lambda^+$ for some $\sigma$ in the classical Weyl group $\Weyl,$
$\gamma \in \Lambda^l,$  suppose
$\lambda$ is on exactly one wall $w_{k,\alpha}$ and
$\lambda+\sigma(\gamma)$ is in the interior of an alcove for which
$w_{k,\alpha}$ is a lower wall.   Then $T_{\lambda+\sigma(\gamma)}$
is a direct summand of $T_\lambda \tensor T_\gamma.$
\end{lemma}
\begin{pf}
By Lemma~\ref{tensor_factors}(b), we must check that
$\lambda+\sigma(\gamma)$ occurs as a highest weight in the classical
decomposition of $W^\lambda \tensor W^\gamma,$ and that nothing
greater in its
$\Weyl_l$ orbit occurs as a highest weight in  classical $W^{\lambda'}
\tensor W^\gamma$ with $\lambda'$ $\Weyl_l$-equivalent to and less
than $\lambda$ (nothing is  $\Weyl_l$-equivalent to and less
than $\gamma$ because it is in the Weyl alcove).

For the first point, consider the classical Racah formula, Equation
(\ref{Racah}). Note that the result is true unless
$\lambda+\sigma(\gamma)$ is the result of the translated action of a
nontrivial element 
$\tau$ of
the classical Weyl group on $\lambda+\mu$ for some $\mu$ that occurs
as a weight of $W^\gamma.$   It is easy to see that if $\lambda,
\lambda'$ are in the Weyl chamber then $\tau\cdot \lambda'$ is
strictly further
from $\lambda$ than $\lambda'$ for any $\tau \in \Weyl,$ so the length
of $\mu$ must be strictly greater than the length of $\sigma(\gamma),$
which is not possible if $\mu$ is a weight of $W^\gamma.$

Essentially the same argument applies for the second point.  Since
$\lambda,\lambda+\sigma(\gamma)$ are in the same alcove,  any
$\lambda',\mu'$ in the translated $\Weyl_l$ orbits respectively of
$\lambda$ and 
$\lambda+\sigma(\gamma)$  must be at least as far
away from each other as $\lambda+\sigma(\gamma)$ and $\lambda$ are, with
equality only when $\mu',\lambda'$ are in the same alcove.  But if $\mu'$ is a
weight in classical $W^{\lambda'} \tensor W^\gamma,$ it must be at
most $||\gamma||$ away from $\lambda',$ with that distance only
achieved if $\mu'-\lambda'=\sigma'(\gamma)$ for some $\sigma' \in
\Weyl.$ Thus if $\lambda'<\lambda,$  $\lambda',\mu'$ are in the orbits
of $\lambda$ and $\lambda + \sigma(\lambda),$ and $\mu'$ is in
$W^\lambda \tensor W^\gamma,$ then $\mu'$ is in the same alcove as
$\lambda',$ so there is a single $\tau \in \Weyl_l$ such that
$\tau\cdot \lambda=\lambda'$ and $\tau \cdot
(\lambda+\sigma(\lambda))=\mu'.$ If $\mu'>\lambda+\sigma(\gamma),$
then $\lambda' \geq \lambda,$ so we must have $\lambda=\lambda'.$  In
this case $\tau=\sigma_{n,\alpha},$ the reflection about the wall on
which $\lambda$ lies.  since $\lambda$ is on a lower wall this would
make $\mu'\leq\lambda+\sigma(\gamma).$  Thus by contradiction the
result is proven.
\end{pf}

\begin{pf}[Of Theorem~\ref{negligible}] In light of
Lemma~\ref{wall_factor} and Lemma~\ref{walls}, it suffices to find for
each alcove  other than
$C_l$ with nonempty intersection with $\Lambda^+$  a
dominant weight $\lambda$ on the interior of a lower wall of that
alcove.  Then
each $\mu$ in the interior of this alcove, since $\mu-\lambda$ is
 Weyl conjugate 
to something in $C_l,$  would have a $T_\mu$ as a summand in some
$T_\lambda \tensor T_\gamma,$ and thus would be negligible. 
This requires that every such alcove have a lower wall whose
intersection with the weight lattice consists of dominant weights, and
that on the interior of each wall of the fundamental
alcove there is a weight.

For the first, notice that every wall of an alcove is either a wall of
the principal chamber for the translated action of the classical Weyl
group or is transverse to it, so every wall of every alcove either
contains no dominant weights or all the weights in its interior are
dominant. If the alcove intersects $\Lambda^+$ all of its walls that do
not intersect $\Lambda^+$ must be part of the walls of the
 chamber.  If all the lower walls of an alcove are walls of
the chamber,  the alcove clearly must be $C_l.$

For the second, if the the wall is $w_{0,\alpha_i},$ one can readily
check that $-\lambda_i$ lies in the interior of the wall (under the
restriction on $l$).  If the wall
is $w_{1,\theta}$ ($D|l'$ case), notice there is always a fundamental weight
$\lambda_i$ such that $\bracket{\lambda_i,\check{\theta}}=1$ (Check
\cite{Humphreys72}[p. 66]), so $(l'-\check{h}/D)\lambda_i$ lies on
$w_{1,\theta}.$ Since it is a dominant weight it lies on no other
walls.   If the wall is $w_{1,\phi}$ ($D\!\! \not{|}l'$ case),  we can find $\lambda_i$ such
that $\bracket{\lambda_i,\phi}=1$ for $B_n$ and $C_n,$ and therefore
$(l'-h)\lambda_i$ will do the trick. There remains only $G_2$ and
$F_4$ to consider.

For $G_2,$ we check that $\bracket{\lambda_1,\phi}=2,$
$\bracket{\lambda_2,\phi}=3.$ Now every integer greater than $1$ can
be written as a nonnegative integer combination of $2$ and $3,$ and every
number greater than $6$ can be written so with neither coefficient
equal to zero.  Thus if $l>6=h$ there exists a positive integer
combination of $\lambda_1$ and $\lambda_2$ whose inner product with
$\phi$ is $l.$  Thus this integer combination minus $\rho$ lies on
$w_{l,\phi}$ and no other wall.

For $F_4,$ we check that $\bracket{\lambda_i,\phi}$ gives $2,4,3,2$
for $i=1\ldots4.$  Again if $l>12-h,$ then $l$ can be written as a
positive linear combination of these four numbers, and thus the same
combination of $\lambda_1$ through $\lambda_4$ gives a weight on the
interior of $w_{l,\phi}.$   
\end{pf}

\section{The Quantum Racah Formula}

For classical Lie algebras or generic $q$ write
\begin{equation} \label{tensor_formula} W^\lambda \tensor W^\gamma
\iso \bigoplus_{\mu \in \Lambda^+} N_{\lambda,\gamma}^\mu W^\mu
\end{equation}
where $N_{\lambda,\gamma}^\mu$ are nonnegative integers representing
multiplicities. 

  For $\q$ an $l$th root of unity, if 
$\lambda,\gamma \in \Lambda^l,$ then  
$$W^\lambda \tensor W^\gamma \iso \bigoplus_{\mu \in \Lambda^l}
M_{\lambda,\gamma}^\mu W^\mu \oplus N$$
where $N$ is a negligible tilting module and each $M_{\lambda,\gamma}^\mu$ is
a nonnegative integer representing multiplicity (because $W^\lambda,$
$W^\gamma,$ and $W^\mu$ are all titling modules). We define the
\emph{truncated tensor product} $\trunc$ on direct sums of Weyl modules in
$\Lambda^l$ by extending the following to direct sums:
\begin{equation}\label{truncated_tensor}W^\lambda \trunc W^\gamma=
W^\lambda \tensor W^\gamma/N \iso 
\bigoplus_{\mu \in \Lambda^l} 
M_{\lambda,\gamma}^\mu W^\mu.\end{equation} 
We will see in the next section that this gives a monoidal structure
on the category of such modules.
 
\begin{proposition}
\begin{equation}
\label{quantum_tensor_decomp}
M_{\lambda,\gamma}^\mu=\sum_{\mu \in \Lambda^l}   \sum_{\begin{matrix}
\sigma \in 
\Weyl_l\\ \sigma\cdot \mu \in \Lambda^+\end{matrix}}(-1)^\sigma
N_{\lambda,\gamma}^{\sigma\cdot \mu}.
\end{equation}
\end{proposition}

\begin{pf}
Over $\A,$ 
$$\qtr_{W^\lambda \tensor W^\gamma}=\sum_\mu N_{\lambda,\gamma}^\mu
\qtr_{W^\mu}$$
so in particular the same holds over $\QQ[\s].$  As a functional on
the center this is equal to
$$\sum_{\mu \in \Lambda^l}   \sum_{\begin{matrix} \sigma \in
\Weyl_l\\ \sigma\cdot \mu \in \Lambda^+\end{matrix}}(-1)^\sigma
N_{\lambda,\gamma}^{\sigma\cdot \mu}.$$
On the other hand as a functional on the center 
$$\qtr_{W^\lambda \tensor W^\gamma}=\qtr_{W^\lambda \trunc
W^\gamma}=\sum_{\mu \in \Lambda^l} M_{\lambda,\gamma}^\mu \qtr_\mu.$$
Since $\{\qtr_\mu\}_{\mu \in \Lambda^l}$ are linearly independent as
functionals on the center (Corollary~\ref{independent}), the result follows.   
\end{pf}
\begin{corollary}[Quantum Racah Formula]
\begin{equation}
\label{quantum_racah}M_{\lambda,\gamma}^\mu=\sum_{\sigma \in \Weyl_l}
(-1)^\sigma \dim(W^\lambda(\sigma\cdot\mu -\gamma)) 
\end{equation}
where $W^\lambda(\gamma)$ is the subspace of $W^\lambda$ of weight
$\gamma.$
\end{corollary}  
\begin{pf}
This result relies on the \emph{classical Racah formula,} which says that 
\begin{equation}\label{Racah}N_{\lambda,\gamma}^\mu=\sum_{\tau \in
\Weyl} (-1)^\tau \dim(W^\lambda(\tau\cdot\mu -\gamma)). 
\end{equation}
  Recall that $\Lambda^+$ is a fundamental domain for the action of
  $\Weyl$ and that only the identity fixes it.  Suppose $\sigma \in
  \Weyl_l$ takes the principal Weyl alcove $C_l$ to some domain $C.$
  There is a unique element $\tau \in \Weyl$ such that
  $\tau^{-1}$ of $C$ intersects $\Lambda^+.$ Thus
  $\tau^{-1}\sigma$ takes $\Lambda^l$ to some fundamental
  domain in $\Lambda^+.$  We conclude that every element of the
  affine Weyl group can be written uniquely as $\tau
  \eta,$ where $\tau \in \Weyl$ and $\eta[\Lambda^l]\subset \Lambda^+.$ Thus
\begin{multline*}
M_{\lambda,\gamma}^\mu=\sum_{\eta[\Lambda^l]\subset \Lambda^+} (-1)^\eta
N_{\lambda,\gamma}^{\eta\cdot\mu} \\= \sum_{\eta[\Lambda^l]\subset \Lambda^+} (-1)^\eta
\sum_{\tau \in \Weyl} (-1)^\tau
\dim(W^\lambda(\tau\eta\cdot \mu-\gamma) \\=\sum_{\sigma \in
\Weyl_l} (-1)^\sigma \dim(W^\lambda(\sigma\cdot \mu -\gamma).
\end{multline*}
\end{pf}

\begin{remark}
The quantum Racah formula (\ref{quantum_racah}), just like the
classical version, admits a beautiful concrete algorithm for the
computation of $M_{\lambda,\gamma}^\mu$ in rank $2$ which illustrates its geometric flavor.  Draw
the weight lattice.  Cover this with a piece of tracing paper and mark
off next to each weight the dimension of the corresponding weight
space of $V^\gamma.$  Now slide the tracing paper so that what
initially covered the the $0$
weight space now lies over the weight space $\lambda.$  Fold the tracing
paper along the walls of $C_l,$ and 
continue folding until the paper fits within it.  For each weight
$\mu$ add up all the numbers that now lie over the point $\mu,$
subtracting those numbers that appear in reverse (write with a
sufficiently seriphed font that you can distinguish them!).  This sum
is $M_{\lambda,\gamma}^\mu.$ 
\end{remark}

\begin{corollary}\label{symmetry}
Let $\iota$ be an isometry of $C_l$ which preserves weights (hence
also an isometry of $\Lambda^l$) given by a translation composed with
the translated action of a Weyl group element.  Then for all
$\lambda,\gamma, \mu
\in \Lambda^l,$ 
$$M_{\lambda,\gamma}^\mu=M_{\lambda,\iota(\gamma)}^{\iota(\mu)}.$$
\end{corollary}

\begin{pf}
Write $\iota(\gamma)=\sigma \cdot \gamma + t,$ where $\sigma \in
\Weyl$ and $t$ is a weight.  Notice
$\sigma(\lambda-\gamma)=\iota(\lambda)-\iota(\gamma).$  Also, since
$\iota(-\rho)=t+\rho$ is a vertex of $C_l,$
$\bracket{t,\check{\alpha}_i}=l'_i,$ which is to say
$\sigma_{0,\alpha_i}(t)-t = l_i'\alpha_i \in M$ for all simple
$\alpha_i \in \Delta.$  Likewise $\sigma_{1,\theta_0}(t)-t \in M$
(check on translation by $l'\theta_0$) so $\sigma(t)-t\in M$ for all
$\sigma \in \Weyl_l.$ Therefore
$$\iota^{-1} \sigma' \iota \in \Weyl_l \qquad \text{ for all }\sigma'
\in \Weyl_l.$$
Then
\begin{multline*}
M_{\lambda,\gamma}^\mu=\sum_{\sigma' \in \Weyl_l} (-1)^{\sigma'}
\dim(W^\lambda(\sigma'\cdot \mu - \gamma))\\
=\sum_{\sigma' \in \Weyl_l} (-1)^{\sigma'}
\dim(W^\lambda(\iota(\sigma'\cdot \mu) - \iota(\gamma)))\\
=\sum_{\sigma' \in \Weyl_l} (-1)^{\sigma'}
\dim(W^\lambda(\sigma''\cdot\iota(\mu) - \iota(\gamma)))\\
=M_{\lambda,\iota(\gamma)}^{\iota(\mu)}$$
\end{multline*}
where $\sigma''$ is the conjugate of $\sigma'$ by $\iota,$ which of
course has the same orientation as $\sigma'.$
\end{pf}

\section{Ribbon Categories and Modularity}

\begin{theorem}\label{ribbon_category} The category of 
modules of $U_{\A'}^\dagger(\lieg)$ which are free and
finite-dimensional over $\A'$ is an Abelian ribbon category enriched over
$\A'.$  In particular the invariant of a link with components labeled
by finite-rank modules is a polynomial in $s$ with integer coefficients.
\end{theorem} 

\begin{pf}
Kassel in \cite{Kassel95}[XI-XIV] defines ribbon categories and proves
that the category of finite-dimensional modules of a ribbon Hopf
algebra over a field forms a ribbon category.  One can easily check
that the proof goes through unchanged for topological Hopf algebras,
and for Hopf algebras over a p.i.d., provided we restrict to free
modules (the freeness is required to define a map from the trivial
module to the tensor product of a module with its dual).  The last
remark is simply the observation that the trivial object in the
category is the module $\A'$ with the counit as action.
\end{pf}

\begin{remark} Le in \cite{Le00a} has proven a stronger result than
the last sentence.  The link invariant (when all labels are Weyl
modules) consists of a term depending
only on the linking matrix and labels times an integer polynomial in $q^2$ and
$q^{-2}.$  In fact with care a similar result holds for the entire
ribbon category.
\end{remark}

\begin{theorem} \label{q_ribbon_category} The category of all
finite-dimensional 
$U_\s^\dagger(\lieg)$-modules is an Abelian ribbon category enriched
over $\QQ[\s],$  and the
full subcategory of tilting modules is a
ribbon subcategory.  The invariant $I(L)$ of any labeled
link $L$ with a 
component labeled by $W^\lambda$ is related to that of $L'$, the same
link with that component labeled by $W^{\sigma\cdot \lambda}$ for $\sigma
\in \Weyl_l,$ by 
$$\label{invariant_relation}I(L')=(-1)^\sigma I(L).$$
\end{theorem}

\begin{pf}
Again the first sentence is a corollary of Kassel's proof. A
subcategory of a ribbon category is ribbon so long as it is 
 closed under tensor product and left duals, and this is the content of
Corollary \ref{tilting}  (together with the obvious fact that the set
of tilting modules is closed under taking duals).

It is shown for example by Kauffman and Radford \cite{KR95b} that to
a one-tangle with components labeled by invariant functionals on a
ribbon Hopf algebra there is associated an element of the center such
that, if all the functionals are quantum traces of modules $W_i$
and the quantum trace of some module $V$ is applied to the action of this
central element on $V,$ the result is the ordinary link
invariant of the link labeled by the given modules, with the open
component closed and labeled by $V.$ In other words the
invariant is $\qdim(V)\chi_V(z),$ where $V$ is the label of the open
component and $z$ is the element of the center.  The result follows
from this fact together with Corollary~\ref{qtr_cor}.
\end{pf}

\begin{theorem} \label{quotient} The category of tilting modules has a
full ribbon functor to a semisimple ribbon category $\Cat$ whose
nonisomorphic simple
objects are the image of the tilting modules with highest weight in
$\Lambda^l.$ The invariant of a ribbon spin network with labels in the
original category is the same as the invariant of the same network
labeled with the functorial image of those labels.
\end{theorem}

\begin{proof}
We use the quotient construction of Mac Lane
\cite{MacLane71}[II.8]. Specifically, if $f,g \in \Hom(V,W),$ we say that
$f\sim g$ if, for all $h \in \Hom(W,V),$ $\qtr_V(hf)=\qtr_V(hg).$
Such an equivalence relation defines a functor to a quotient category
$\Cat$ such that 
$f\sim g$ implies $f$ and $g$ are sent by the functor to the same
thing, and $\Cat$ is universal for this property.  

It is clear that the image of any negligible tilting modules is null,
i.e. isomorphic to the null module $\{0\}.$  since $\qdim(\lambda)\neq
0$ for all $\lambda \in \Lambda^l,$ each such module is mapped to a
non-null object and thus $\Cat$ is semisimple with these as
simple objects.  

Because $f \sim g$ implies $f \tensor h \sim g \tensor h$ and $h
\tensor f \sim h \tensor g$ for all $h$, $\Cat$ inherits a 
tensor product structure making the functor a tensor functor.  The
image of the braiding morphisms, duality morphisms and twist morphisms
are braiding, duality and twist morphisms for $\Cat.$
   
\end{proof}
\begin{remark}
Since each tilting module can be written uniquely as a direct sum $M \oplus
N,$ where $M$ is isomorphic to a direct sum of modules in $\Lambda^l$
and $N$ is negligible, we can define an isomorphic functor from 
$\Cat$  to the full subcategory of the tilting module category
consisting of modules
isomorphic to a direct sum of modules in $\Lambda^l$. Unfortunately this
functor does not preserve the tensor product.  However this is a
tensor isomorphism if the range category uses the truncated tensor
product $\trunc$ for its monoidal structure.  This proves the
truncated tensor product is a monoidal structure, and that in
particular each $M_{\lambda,\gamma}^\mu$ is nonnegative.
\end{remark}

To address the modularity of $\Cat$ we must come to terms with the
relationship between $\Weyl_l,$ the symmetries of the characters of
Weyl modules, and $\Weyl^\dagger,$
the symmetries of the $S$-matrix.  Actually, since we are only
interested in the $S$-matrix applied to weights, the relevant group is the
subgroup $\Weyl_l^\Lambda \subset \Weyl^\dagger$ which map weights to weights.

\begin{proposition}
\item $\frac{l}{2} \check{\Lambda}_r \cap \Lambda$ is the lattice
      generated by $M,$ together with the vectors $l'\lambda_i,$ where
$\lambda_i$ is as in Table \ref{table}, using the conventions in
\cite{Humphreys72}[Ch. 13] for the naming of fundamental weights
(when no $\lambda_i$ is given the lattice is exactly $M$).
Notice each $\lambda_i$ corresponds to an element of the fundamental
group of order $2,$ though not every such occurs on the list.
\end{proposition} 

\begin{table} \caption{Instances of extra symmetry in the Weyl
alcove}\label{table} 
\begin{tabular}{|r||c|c|c|c|c|} 
\hline
 & $A_n$ & $B_n$ & $C_n$ & $D_n$ & $E_7$\\
\hline
$l$ odd & $\lambda_{\frac{n+1}{2}},$ if $2\!\!\not|\,n$ & 
$\lambda_n$ if $2|n$
 &
$\lambda_1$ &
\begin{tabular}{c} 
$\lambda_1$ if $2\!\!\not|\,n$\\
$\lambda_1,\lambda_{n-1},\lambda_n$ if 
$2|n$
\end{tabular} & 
$\lambda_7$\\
\hline 
$l$ even & & $\lambda_n$ if $D\!\!\not|\,l',\, 2\!\!\not|\,n$ &
$\lambda_1$  if $D\!\!\not|\,l'$ & &  \\
\hline
\end{tabular}
\end{table}
\begin{pf}
Suppose first that $l$ is odd, and suppose $\check{\alpha} \in
      \check{\Lambda}_r$ such that $l\check{\alpha}/2 \in \Lambda$ but
      is not in $M.$  Since $D\check{\alpha} \in \Lambda_r$ we can
      subtract any multiple of it from  $l\check{\alpha}/2$ and it
      will still be in the weight lattice and still not be in $M.$  As
      long as $D \neq 3$ this 
      means $\check{\alpha}/2 \in \Lambda,$ is half a coroot, and $l$
      times it is not in $M.$  If $D|l$ this means it is a weight
      $\lambda$ such that $2\lambda$ is a coroot but $l\lambda$ is not
      a coroot.  Since $l$ is odd this is equivalent to saying
      $\lambda$ is not a coroot. Since the set of $\lambda$ with these
      properties  
      form a $\Lambda_r$ equivalence class, it suffices to check which
      minimal $\lambda_i$ are half a coroot but not a coroot
      themselves. If $D$ does not divide $l$ the same reasoning says
      we are looking for minimal $\lambda_i$ such that $2\lambda_i$ is
      a coroot but $l\lambda_i$ is not a root. Again assuming $D$ is
      not $3,$ we see that $l \lambda_i$ is a root if and only if
      $\lambda_i$ is a root. The table lists those $\lambda_i$
      which are half a coroot but not a root.

      the case $D=3,$ $G_2,$ is handled similarly.

      If $l$ is even, the two lattices are equal unless $D$ does not
      divide $l'.$  If $l'\check{\alpha}$ is a weight, then by
      subtracting an appropriate multiple of $D\check{\alpha}$ we see
      that $\check{\alpha}$ is a weight.  Arguing as above this is
      equivalent to finding a minimal weight $\lambda_i$ which is a
      coroot but is not a root.  These are listed above.

\end{pf} 

\begin{proposition}
For each $\lambda_i$ in the table above, there is an element
$\sigma_i,$ of the classical
Weyl group which permutes all the simple roots except $\alpha_i$ and
and sends $\alpha_i$ to $-\theta_0.$  This Weyl group element
composed with the translation operator $l'\lambda_i$ is an isometry of
$C_l$ and $\Lambda^l,$ and any two elements of $\Lambda^l$ which have
proportional entries in the $S$-matrix are related by a product of
such isometries.
\end{proposition}

\begin{pf}
If $D|l',$ the isometries of the standard Weyl alcove are discussed in
\cite{Sawin02a}, and shown to be of the form described in the
proposition, with a set of $\lambda_i$ that includes the entries in
the table.

If $D$ does not divide $l',$ we are necessarily in the case $B_n$ and
$C_n.$ One checks that there is a Weyl group element of the sort
described in the proposition (here $\theta_0=\phi,$ and note
that every root is positive or negative with respect to the resulting
base). It is easy to see the inequalities defining $C_l$ are preserved
by the given transformation, and since the transformation sends
weights to weights it must preserve $\Lambda^l.$  

Two elements of $\Lambda^l$ which have proportional entries in the
$S$-matrix must have linearly dependent characters when applied to the
image of $\Psi$ because
$$S_{\lambda,\gamma}=\qdim(\lambda)\chi_\lambda(\Psi(\qtr_\gamma)).$$
Thus two such elements  must be related by an element of
$\Weyl^\dagger.$  Since such a transformation 
preserves $\Lambda,$ it must 
be a product of an ordinary Weyl transformation and a translation
generated by $M$ and $l'\lambda_i.$  Conjugating by ordinary Weyl
transformations as necessary, we can write this as a product of terms
of the form $T_{l'\lambda_i} \sigma_i,$ times an element of $\Weyl_l,$
where $T_\lambda$ is translation by $\lambda.$   
Since the element of $\Weyl_l$ takes $C_l$ to itself, it must be the
identity.
\end{pf}

\begin{theorem}
The category $\Cat$ is modular and thus gives a $2+1$-dimensional TQFT
except in the cases where $\lambda_i$ is given  in Table \ref{table}. 
 In these
cases the quotient category as described by  Brugui\`eres
\cite{Bruguieres00} is well-defined, and gives either a modular or a
spin modular category, except in the cases (with $l$ odd) $A_n$ with
$n \equiv 1 
\bmod 4,$ $B_n$ with $n \equiv 2 
\bmod 4,$ $C_n,$ $D_n$ with $n \equiv 2 
\bmod 4,$ and $E_7$ and (with $l$ even) $B_n$ $l'$ odd and  $n \equiv 1 
\bmod 4,$ and $C_n$ with $l'$ odd.
\end{theorem}

\begin{pf}
Brugui\`eres shows that a semisimple ribbon category is modular unless
there is a simple object whose entries in the $S$-matrix are
proportional to those of the trivial object.  Thus by the previous
proposition those without the special isometries of the Weyl alcove
are modular.

In order for there to be a well-defined minimal quotient which is
modular or spin modular according to Brugui\`eres (the notion of
spin-modular is introduced in \cite{Sawin02b}, but the work of
Brugui\`eres extends directly to it) it suffices to show that all
simple objects whose rows in the $S$-matrix are proportional to that
of the
trivial object have quantum dimension $1,$ are \emph{transparent}
(i.e. the square of the $R$-matrix acts as the identity on the tensor
product of this module with any other), and the set of them forms a
group under tensor product.

By the previous proposition these objects are those of the form
$\iota(0)$ where $\iota$ is one of the isometries associated to the
$\lambda_i$ of Table \ref{table}.  Since these isometries are the
composition of a Weyl group element and a translation, Corollary
\ref{symmetry} tells us $\iota(\lambda)=\lambda \trunc \iota(0)$ (for
simplicity we refer to the object indexed by $\lambda \in \Lambda^l$
as $\lambda$) so these objects form a group under tensor product (in fact the
group $\Weyl_l^\Lambda/\Weyl_l$).

The square of the $R$ matrix as a map from $W^\lambda \trunc
W^{\iota(0)}\iso W^{\iota(\lambda)}$ is a multiple of the identity,
since this module is simple.  Its quantum trace is $S_{\lambda,
\iota(0)}=\qdim(\lambda)\qdim(\iota(0))$ and thus it is the
identity. So $\iota(0)$ is transparent.

Finally a check of the quantum dimension formula
(\ref{Weyl_character}) in each case in the table shows that the
dimensions of the modules of weight $\iota(0)=l' \lambda_i - \rho
+\sigma(\rho)$ is $-1$ in exactly the case given in the theorem, and
$1$ in the other cases.

\end{pf}
\begin{remark}
Presumably, some topological information can still be gleaned in the
cases not covered by the theorem (which include $U_q(\mathfrak{su}_2)$
at odd roots of unity)
\end{remark}

\begin{remark}
The case where $l$ is a multiple of $2D$ corresponds to Chern-Simons
theory at integer level, and is the case of the most physical
interest.  It also is the simplest, and indeed most of the technical
details in this paper deal with the other case.  The case where $l$ is
not a multiple of $2D$ corresponds to certain fractional levels, and
it is not clear from the physical interpretation why we should expect
modularity at these levels.  The fact that it happens at all, as well
as the new behavior these fractional levels exhibit (such as  being
defined over the alcove of the \emph{dual}  affine Weyl group,
and the need in certain case to quotient to achieve modularity), are
phenomena that beg an interpretation in terms of Chern-Simons theory.
\end{remark}
\def\cprime{$'$}

\end{document}